\documentclass[12pt]{article}
\usepackage{amsfonts,amscd,a4}
\usepackage{color}
\usepackage{graphicx}
\usepackage{theorem}
\newtheorem{theo}{Theorem}
\newtheorem{prop}[theo]{Proposition}
\newtheorem{lemma}[theo]{Lemma}
\newtheorem{coro}[theo]{Corollary}

{\theorembodyfont{\rm}

}

\newcommand{\cA}{{\mathcal A}}
\newcommand{\cB}{{\mathcal B}}
\newcommand{\cC}{{\mathcal C}}

\newcommand{\cF}{{\mathcal F}}

\newcommand{\cI}{{\mathcal I}}
\newcommand{\cJ}{{\mathcal J}}

\newcommand{\cP}{{\mathcal P}}

\newcommand{\eH}{{\sf H}}

\newcommand{\eT}{{\sf T}}

\newcommand{\sC}{{\mathbb C}}

\newcommand{\sI}{{\mathbb I}}

\newcommand{\sN}{{\mathbb N}}

\newcommand{\sR}{{\mathbb R}}

\newcommand{\sT}{{\mathbb T}}

\newcommand{\sZ}{{\mathbb Z}}
\newcommand{\qed}{\rule{1ex}{1ex}}

\newcommand{\clos}{\mbox{\rm clos} \,}
\newcommand{\coker}{\mbox{\rm coker} \,}

\newcommand{\diag}{\mbox{\rm diag} \,}

\newcommand{\el}{\mbox{\rm el}}

\newcommand{\ext}{\mbox{\rm ext}}

\newcommand{\im}{\mbox{\rm im} \,}
\newcommand{\ind}{\mbox{\rm ind} \,}

\newcommand{\smb}{\mbox{\rm smb}_p \,}

\newcommand{\wind}{\mbox{\rm wind} \,}

\begin{document}
\title{A handy formula for the Fredholm index of Toeplitz plus Hankel operators}
\author{Steffen Roch and Bernd Silbermann}
\date{Dedicated to the memory of Israel Gohberg}
\maketitle
\begin{abstract}
We consider Toeplitz and Hankel operators with piecewise continuous generating
functions on $l^p$-spaces and the Banach algebra generated by them. The goal of
this paper is to provide a transparent symbol calculus for the Fredholm property
and a handy formula for the Fredholm index for operators in this algebra.
\end{abstract}
{\bf Keywords:} Toeplitz plus Hankel operators, Fredholm index \\
{\bf 2010 AMS-MSC:} Primary: 47B35, secondary: 47B48
\section{Introduction} 
Throughout this paper, let $1 < p < \infty$. For a non-empty subset $\sI$ of the
set $\sZ$ of the integers, let $l^p (\sI)$ denote the complex Banach space of
all sequences $x = (x_n)_{n \in \sI}$ of complex numbers with norm $\|x\|_p =
(\sum_{n \in \sI} |x_n|^p)^{1/p} < \infty$. We consider $l^p(\sI)$ as a closed
subspace of $l^p(\sZ)$ in the natural way and write $P_\sI$ for the canonical
projection from $l^p(\sZ)$ onto $l^p(\sI)$. For $\sI = \sZ^+$, the set of the
non-negative integers, we write $l^p$ and $P$ instead of $l^p(\sI)$ and $P_\sI$,
respectively. By $J$ we denote the operator on $l^p(\sZ)$ acting by $(Jx)_n :=
x_{-n-1}$, and we set $Q := I - P$.

For every Banach space $X$, let $L(X)$ stand for the Banach algebra of all
bounded linear operators on $X$, and write $K(X)$ for the closed ideal of $L(X)$
of all compact operators. The quotient algebra $L(X)/K(X)$ is known as the
Calkin algebra of $X$. Its importance in this paper stems from the fact that the
invertibility of a coset $A + K(X)$ of an operator $A \in L(X)$ in this
algebra is equivalent to the Fredholm property of $A$, i.e., to the finite
dimensionality of the kernel $\ker A = \{x \in X: Ax = 0\}$ and the cokernel
$\coker A = X/\im A$ of $A$, with $\im A = \{Ax : x \in X\}$ referring to the
range of $A$. If $A$ is a Fredholm operator then the difference $\ind A :=
\dim \ker A - \dim \coker A$ is known as the Fredholm index of $A$.

Our goal is a criterion for the Fredholm property and a formula for the Fredholm
index for operators in the smallest closed subalgebra of $L(l^p)$ which contains
all Toeplitz and Hankel operators with piecewise continuous generating function.
The precise definition is as follows. Let $\sT$ be the complex unit circle. For
each function $a \in L^\infty(\sT)$, let $(a_k)_{k \in \sZ}$ denote the sequence
of its Fourier coefficients,
\[
a_k := \frac{1}{2 \pi} \int_0^{2 \pi} a(e^{i\theta}) e^{-ik\theta} \, d\theta.
\]
The {\em Laurent operator} $L(a)$ associated with $a \in L^\infty(\sT)$ acts on
the space $l^0(\sZ)$ of all finitely supported sequences on $\sZ$ by
$(L(a) x)_k := \sum_{m \in \sZ} a_{k - m} x_m$. (For every $k \in \sZ$, there
are only finitely many non-vanishing summands in this sum.) We say that $a$ is a
multiplier on $l^p(\sZ)$ if $L(a) x \in l^p(\sZ)$ for every $x \in l^0(\sZ)$ and
if
\[
\|L(a)\| := \sup \{ \|L(a)x\|_p : x \in l^0(\sZ), \, \|x\|_p = 1\}
\]
is finite. In this case, $L(a)$ extends to a bounded linear operator on
$l^p(\sZ)$ which we denote by $L(a)$ again. The set $M^p$ of all multipliers on
$l^p(\sZ)$ is a Banach algebra under the norm $\|a\|_{M_p} := \|L(a)\|$. We let
$M^{\langle p \rangle}$ stand for $M^2$ if $p=2$ and for the set of all $a \in
L^\infty(\sT)$ which belong to $M^r$ for all $r$ in a certain open neighborhood
of $p$ if $p \neq 2$.

It is well known that $M^2 = L^\infty(\sT)$. Moreover, every function $a$ with
bounded total variation $\mbox{Var} (a)$ is in $M^p$ for every $p$, and the
Stechkin inequality
\[
\|a\|_{M_p} \le c_p (\|a\|_\infty + \mbox{Var} (a))
\]
holds with a constant $c_p$ independent of $a$. In particular, every
trigonometric polynomial and every piecewise constant function on $\sT$ are
multipliers for every $p$. We denote the closure in $M^p$ of the algebra $\cP$
of all trigonometric polynomials and of the algebra $P\sC$ of all piecewise
constant functions by $C_p$ and $PC_p$, respectively. Thus, $C_p$ and $PC_p$ are
closed subalgebras of $M^p$ for every $p$. Note that $C_2$ is just the algebra
$C(\sT)$ of all continuous functions on $\sT$, and $PC_2$ is the algebra
$PC(\sT)$ of all piecewise continuous functions on $\sT$. It is well known that
$C_p \subseteq C(\sT)$ and $C_p \subseteq PC_p \subseteq PC(\sT)$ for every $p$.
In particular, every multiplier $a \in PC_p$ possesses one-sided limits
at every point $t \in \sT$ (see \cite{BS1} for these and further properties of
multipliers). For definiteness, we agree that $\sT$ is oriented
counter-clockwise, and we denote the one-sided limit of $a$ at $t$ when
approaching $t$ from below (from above) by $a(t^-)$ (by $a(t^+)$).

Let $a \in M^p$. The operators $T(a) := P L(a) P$ and $H(a) := P L(a) Q J$,
thought of as acting on $\im P = l^p$ are called the Toeplitz and Hankel
operator with generating function $a$, respectively. It is well known
that $\|T(a)\| = \|a\|_{M_p}$ and $\|H(a)\| \le \|a\|_{M_p}$
for every multiplier $a \in M_p$.

For a subalgebra $A$ of $M^p$, we let $\eT(A)$ and $\eT\eH (A)$ stand for the
smallest closed subalgebra of $L(l^p)$ which contains all operators $T(a)$ with
$a \in A$ and all operators $T(a) + H(b)$ with $a, \, b \in A$, respectively. We
will be mainly concerned with the algebras $C_p$, $PC_p$, and with their
intersections with $M ^{\langle p \rangle}$, in place of $A$. Now we can state
the goal of the paper more precisely: we will state a criterion for the Fredholm
property of operators in $\eT\eH (PC_p)$ and derive a formula for the Fredholm
index of operators $T(a) + H(b)$ with $a, \, b \in PC_p$.

The study of the Fredholm property of operators in $\eT\eH (PC_p)$ has a long
and involved history. We are going to mention only some of its main stages.

The Fredholm properties of operators in the algebra $\eT (PC_p)$ are well
understood thanks to the work of I. Gohberg/N. Krupnik and R. Duduchava; see
\cite{BS1} and the literature cited there. We will need these results later on;
therefore we recall them in Section \ref{s2}. Different approaches to these
algebras were developed in \cite{BS1} and \cite{HRS1}; our presentation will be
mainly based on the latter.

The structure of the algebras $\eT\eH (PC_p)$ is much more involved than that of
$\eT (PC_p)$. For instance, the Calkin image $\eT^\pi (PC) := \eT (PC)/K(l^2)$
of $\eT(PC)$ is a commutative algebra, whereas that one of $\eT\eH (PC)$ is not.
The Calkin image of $\eT\eH (PC)$ was first described by Power \cite{Po1}. An
alternative approach was developed by one of the authors in \cite{Sil1}, where
it was shown that the algebra $\eT \eH^\pi (PC) := \eT \eH (PC)/K(l^2)$
possesses a matrix-valued Fredholm symbol. In the present paper, we take up the
approach from \cite{Sil1} in order to study the Fredholm properties of operators
in $\eT\eH (PC_p)$ for $p \neq 2$.

It should be mentioned that the algebras $\eT\eH (PC_p)$ have close relatives
which live on other spaces than $l^p$, such as the Hardy spaces $H^p(\sR)$ and
the Lebesgue spaces $L^p(\sR^+)$. The corresponding algebras were examined (with
different methods) in the report \cite{RoS1}, see also the recent monograph
\cite{RSS1}. Despite these fairly complete results for the Fredholm property, a
general, transparent and satisfying formula for the Fredholm index of operators
in $\eT\eH (PC_p)$ (or on related algebras) was not available until now.
Among the particular results which hold under special assumptions we would like
to emphasize the following. In \cite{KaS1}, there is derived an index formula
for operators of the form $\lambda I + H$ where $\lambda \in \sC$ and $H$ is a
Hankel operator on $H^p(\sR)$. Already earlier, some classes of Wiener-Hopf plus
Hankel operators were studied in connection with diffraction problems; see
\cite{LMT1,MST1}. Note also that the (very hard) invertibility problem for
Toeplitz plus Hankel operators is treated in \cite{BaE1,Ehr1}.

Finally we would like to mention that algebras like $\eT\eH (PC_p)$ can also be
viewed of as subalgebras of algebras generated by convolution-type operators
and Carleman shifts changing the orientation. First results in that direction
were presented in \cite{GK8,GK9} where, in particular, a matrix-valued Fredholm
symbol was constructed.

The goal of the present paper is to provide a transparent symbol calculus for
the Fredholm property as well as a handy formula for the Fredholm index for
operators in the algebra $\eT\eH (PC_p)$. The techniques developed and used in
this paper also allow to handle the corresponding questions for the related
algebras on the spaces $H^p(\sR)$ and $L^p(\sR^+)$.
\section{The Fredholm property} \label{s2}
In what follows, we fix $p \in (1, \, \infty)$ and consider all operators as
acting on $l^p$ unless stated otherwise.

As already mentioned, we start with recalling the basic results of the Fredholm
theory of operators in the algebra $\eT (PC_p)$, which are due Gohberg/Krupnik
and Duduchava. The functions $f_{\pm 1} (t) := t^{\pm 1}$ are multipliers for
every $p$. It is easy to check that the algebra generated by the Toeplitz
operators $T(f_{\pm 1})$ contains a dense subalgebra of $K(l^p)$. Thus, the
ideal $K(l^p)$ is contained in $\eT(C_p)$, hence also in $\eT(PC_p)$, and it
makes sense to consider the quotient algebra $\eT(PC_p)/K(l^p)$. Clearly, if $A
\in \eT(PC_p)$ and if the coset $A + L(l^p)$ is invertible in
$\eT(PC_p)/K(l^p)$, then it is also invertible in the Calkin algebra
$L(l^p)/K(l^p)$, hence $A$ is a Fredholm operator. The more interesting question
is if the converse holds, i.e., if the invertibility of $A + L(l^p)$ in the
Calkin algebra implies the invertibility of $A + K(l^p)$ in $\eT(PC_p)/K(l^p)$.
If this implication holds for every $A \in \eT(PC_p)$, one says that
$\eT(PC_p)/K(l^p)$ is inverse closed in $L(l^p)/K(l^p)$.

Let $\overline{\sR}$ denote the two-point compactification of the real line by
the points $\pm \infty$ (thus $\overline{\sR}$ is homeomorphic to a closed
interval) and let the function $\mu_p : \overline{\sR} \to \sC$ be defined by
\[
\mu_p(\lambda) := (1 + \coth (\pi (\lambda + i/p)))/2
\]
if $\lambda \in \sR$ and by $\mu_p(- \infty) = 0$ and $\mu_p(+ \infty) = 1$.
Note that when $\lambda$ runs from  $- \infty$ to $\infty$ then $\mu_p(\lambda)$
runs along a circular arc in $\sC$ which joins 0 to 1 and passes through the
point $(1 - i \cot(\pi/p))/2$. An easy calculation gives $\mu_p(-\lambda) =
1- \mu_q(\lambda)$, where $1/p + 1/q = 1$. Thus, for fixed $t \in \sT$, the
values $\Gamma (T(a) + K(l^p)) (t, \, \lambda)$ defined in the following
theorem run from $a(t-0)$ to $a(t+0)$ along a circular arc when $\lambda$ runs
from $- \infty$ to $\infty$.
\begin{theo} \label{t1}
$(a)$ $\eT(PC_p)/K(l^p)$ is a commutative unital Banach algebra. \\[1mm]
$(b)$ The maximal ideal space of $\eT(PC_p)/K(l^p)$ is homeomorphic with the
cylinder $\sT \times \overline{\sR}$, provided with an exotic (non-Euclidean)
topology. \\[1mm]
$(c)$ The Gelfand transform $\Gamma : \eT(PC_p)/K(l^p) \to C(\sT \times
\overline{\sR})$ of the coset $T(a) + K(l^p)$ with $a \in PC_p$ is
\[
\Gamma (T(a) + K(l^p)) (t, \, \lambda) = a(t-0) (1 - \mu_q(\lambda)) + a(t+0)
\mu_q(\lambda).
\]
$(d)$ $\eT(PC_p)/K(l^p)$ is inverse closed in $L(l^p)/K(l^p)$.
\end{theo}
The topology mentioned in assertion $(b)$ will be explicitly described in
Section \ref{s3}. Note that this topology is independent of $p$. Since the
cosets $T(a) + K(l^p)$ with $a \in PC_p$ generate the algebra
$\eT(PC_p)/K(l^p)$, the Gelfand transform on $\eT(PC_p)/K(l^p)$ is completely
described by assertion $(c)$. Thus, if $A \in \eT(PC_p)$, then the coset $A +
K(l^p)$ is invertible in $\eT(PC_p)/K(l^p)$ if and only if the function
$\Gamma(A + K(l^p))$ does not vanish on $\sT \times \overline{\sR}$. Together
with assertion $(d)$ this shows that $A \in \eT(PC_p)$ is a Fredholm operator if
and only if $\Gamma(A + K(l^p))$ does not vanish on $\sT \times \overline{\sR}$.
It is therefore justified to call the function $\smb A := \Gamma(A + K(l^p))$
the {\em Fredholm symbol} of $A$.

The index of a Fredholm operator in $\eT(PC_p)$ can be determined my means of
its Fredholm symbol. First suppose that $a \in PC_p$ is a piecewise smooth
function with only finitely many jumps. Then the range of the function
\[
\Gamma (T(a) + K(l^p))(t, \, \lambda) = a(t^-) (1-\mu_q(\lambda)) + a(t^+)
\mu_q)(\lambda)
\]
is a closed curve with a natural orientation, which is obtained from the
(essential) range of $a$ by filling in the circular arcs
\[
\cC_q(a(t^-), \, a(t^+)) := \{ a(t^-) (1-\mu_q(\lambda)) + a(t^+)
\mu_q)(\lambda) : \lambda \in \overline{\sR} \}
\]
at every point $t \in \sT$ where $a$ has a jump. (If the function $a$ is
continuous at $t$, then $\cC_q(a(t^-), \, a(t^+))$ reduces to the
singleton $\{a(t)\}$.) If this curve does not pass through the origin, then we
let $\wind \Gamma (T(a) + K(l^p))$ denote its winding number with respect to the
origin, i.e., the integer $1/(2\pi)$ times the growth of the argument of $\Gamma
(T(a) + K(l^p))$ when $t$ moves along $\sT$ in positive (= counter-clockwise)
direction. If this condition is satisfied then $T(a)$ is a Fredholm operator,
and
\[
\ind T(a) = - \wind \Gamma (T(a) + K(l^p))
\]
(see \cite{BS1}, Section 2.73 and Proposition 6.32 for details). Moreover, as
in Section 5.49 of \cite{BS1}, one can extend both the definition of the
winding number and the index identity to the case of an arbitrary Fredholm
operator in $\eT(PC_p)$. More precisely, one has the following.
\begin{prop} \label{pneu1}
Let $A \in \eT(PC_p)$ be a Fredholm operator. Then
\[
\ind A = -\wind \Gamma (A + K(l^p)).
\]
\end{prop}
We would like to emphasize an important point. The algebra $\eT(PC_2)/K(l^2)$ is
a commutative $C^*$-algebra, hence the Gelfand transform is an isometric
$^*$-isomorphism from $\eT(PC_2)/K(l^2)$ onto $C(\sT \times \overline{\sR})$. In
particular, the radical of $\eT(PC_2)/K(l^2)$ is trivial, and the equality
$\mbox{smb}_2 \, A = 0$ for some operator $A \in \eT(PC_2)$ implies that $A$ is
compact. For general $p$ it is not known if the radical of $\eT(PC_p)/K(l^p)$ is
still trivial; it is therefore not known if $\smb A = 0$ implies the compactness
of $A$.

In order to state our results on the Fredholm property of operators in the
Toeplitz+Hankel algebra $\eT\eH(PC_p)/K(l^p)$ we need some notation. Let $\sT_+$
be the set of all points in $\sT$ with non-negative imaginary part and set
$\sT_+^0 := \sT_+ \setminus \{-1, \, 1\}$. Further let the function $\nu_p :
\overline{\sR} \to \sC$ be defined by
\[
\nu_p(\lambda) := (2i \,\sinh (\pi (\lambda + i/p)))^{-1}
\]
if $\lambda \in \sR$ and by $\nu_p(\pm \infty) = 0$. Recall that $1/p + 1/q =
1$.
\begin{theo} \label{t2}
$(a)$ Let $a, \, b \in PC_p$. Then the operator $T(a) + H(b)$ is Fredholm if
and only if the matrix
\begin{eqnarray} \label{e3}
\lefteqn{\smb (T(a) + H(b)) (t, \, \lambda) :=} \\
&& \pmatrix{a(t^+) \mu_q(\lambda) + a(t^-)(1-\mu_q(\lambda)) &
(b(t^+)-b(t^-)) \nu_q(\lambda) \cr
(b(\bar{t}^-)- b(\bar{t}^+)) \nu_q(\lambda) &
a(\bar{t}^-) (1-\mu_q(\lambda)) + a(\bar{t}^+) \mu_q(\lambda)} \nonumber
\end{eqnarray}
is invertible for every $(t, \, \lambda) \in \sT_+^0 \times \overline{\sR}$ and
if the number
\begin{eqnarray} \label{e4}
\lefteqn{\smb (T(a) + H(b)) (t, \, \lambda) :=} \\
&& a(t^+) \mu_q(\lambda) + a(t^-)(1-\mu_q(\lambda)) + it \, (b(t^+)-b(t^-))
\nu_q(\lambda) \nonumber
\end{eqnarray}
is not zero for every $(t, \, \lambda) \in \{\pm 1\} \times \overline{\sR}$.
\\[1mm]
$(b)$ The mapping $\smb$ defined in assertion $(a)$ extends to a continuous
algebra homomorphism from $\eT\eH(PC_p)$ to the algebra $\cF$ of all bounded
functions on $\sT_+ \times \overline{\sR}$ with values in $\sC^{2 \times 2}$ on
$\sT_+^0 \times \overline{\sR}$ and with values in $\sC$ on $\{\pm 1\} \times
\overline{\sR}$. Moreover, there is a constant $M$ such that
\begin{equation} \label{e241111.1}
\| \smb A \| := \sup_{(t, \lambda) \in \sT_+ \times \overline{\sR}} \|\smb A 
(t, \, \lambda)\|_\infty \le M \inf_{K \in K(l^p)} \|A + K\|
\end{equation}
for every operator $A \in \eT\eH(PC_p)$. Here, $\|B\|_\infty$ refers to the
spectral norm of the matrix $B$.\\[1mm]
$(c)$ An operator $A \in \eT\eH(PC_p)$ has the Fredholm property if and only if
the function $\smb A$ is invertible in $\cF$. \\[1mm]
$(d)$ The algebra $\eT\eH(PC_p)/K(l^p)$ is inverse closed in $L(l^p)/K(l^p)$.
\end{theo}
Before going into the details of the proof, we remark two consequences of
Theorem \ref{t2} which will be needed in the next section.
\begin{coro} \label{c2a}
Let $a, \, b \in PC_p$ and $T(a) + H(b)$ a Fredholm operator on $l^p$.
Then \\[1mm]
$(a)$ the function $a$ is invertible in $PC_p$, and \\[1mm]
$(b)$ if $b$ is continuous at $\pm 1$, then $T(a) - H(b)$ is a Fredholm operator
on $l^p$.
\end{coro}
{\bf Proof.} If $T(a) + H(b)$ is a Fredholm operator, then the diagonal
matrices
\[
\smb (T(a) + H(b)) (t, \, \pm \infty) = \diag (a(t^\pm), \, a(\overline{t}^\pm))
\]
are invertible for every $t \in \sT^0_+$ and the numbers $\smb (T(a) + H(b)) (1,
\, \pm \infty) = a(1^\pm)$ and $\smb (T(a) + H(b)) (-1, \, \pm \infty) =
a((-1)^\pm)$ are not zero by assertion $(a)$ of Theorem \ref{t2}. Hence, $a$ is
invertible as an element of $PC$. Since the algebra $PC_p$ is inverse closed in
$PC$ by Proposition 6.28 in \cite{BS1}, assertion $(a)$ follows. The proof of
assertion $(b)$ is also immediate from the form of the symbol described in
Theorem \ref{t2} $(a)$. \hfill \qed \\[3mm]
The remainder of this section is devoted to the proof of Theorem \ref{t2}.
We will need two auxiliary ingredients which we are going to recall first. Let
$\cA$ be a unital Banach algebra. The {\em center} of $\cA$ is the set of all
elements $a \in \cA$ such that $ab=ba$ for all $b \in \cA$. A {\em central}
subalgebra of $\cA$ is a closed subalgebra $\cC$ of the center of $\cA$ which
contains the identity element. Thus, $\cC$ is a commutative Banach algebra with
compact maximal ideal space $M(\cC)$. For each maximal ideal $x$ of $\cC$,
consider the smallest closed two-sided ideal $\cI_x$ of $\cA$ which contains
$x$, and let $\Phi_x$ refer to the canonical homomorphism from $\cA$ onto the
quotient algebra $\cA/\cI_x$.

In contrast to the commutative setting, where $\cC/x \cong \sC$ for all $x \in
M(\cC)$, the quotient algebras $\cA/\cI_x$ will depend on $x \in M(\cC)$ in
general. In particular, it can happen that $\cI_x = \cA$ for certain maximal
ideals $x$. In this case we {\em define} that $\Phi_x(a)$ is invertible
in $\cA/\cI_x$ for every $a \in \cA$.
\begin{theo}[Allan's local principle] \label{t5}
Let $\cC$ be a central subalgebra of the unital Banach algebra $\cA$. Then an
element $a \in \cA$ is invertible if and only if the cosets $\Phi_x(a)$ are
invertible in $\cA/\cI_x$ for each $x \in M(\cC)$.
\end{theo}
Here is the second ingredient. Recall that an idempotent is an element $p$ of an
algebra such that $p^2 = p$.
\begin{theo}[Two idempotents theorem] \label{t6}
Let $\cA$ be a Banach algebra with identity element $e$, let $p$ and $q$ be
idempotents in $\cA$, and let $\cB$ denote the smallest closed subalgebra of
$\cA$ which contains $p, \, q$ and $e$. Suppose that $0$ and $1$ belong to the
spectrum $\sigma_\cB(pqp)$ of $pqp$ in $\cB$ and that $0$ and $1$ are cluster
points of that spectrum. Then \\[1mm]
$(a)$ for each point $x \in \sigma_\cB(pqp)$, there is a continuous algebra
homomorphism $\Phi_x : \cB \to \sC^{2 \times 2}$ which acts at the generators of
$\cB$ by
\[
\Phi_x(e) = \pmatrix{1 & 0 \cr 0 & 1}, \quad \Phi_x(p) = \pmatrix{1 & 0 \cr 0 &
0}, \quad \Phi_x(q) = \pmatrix{x & \sqrt{x(1-x)} \cr \sqrt{x(1-x)} & 1-x}
\]
where $\sqrt{x(1-x)}$ denotes any complex number with $(\sqrt{x(1-x)})^2 =
x(1-x)$. \\[1mm]
$(b)$ an element $a \in \cB$ is invertible in $\cB$ if and only if the matrices
$\Phi_x(a)$ are invertible for every $x \in \sigma_\cB(pqp)$. \\[1mm]
$(c)$ if $\sigma_\cB(pqp) = \sigma_\cA(pqp)$, then $\cB$ is inverse closed in
$\cA$.
\end{theo}
We proceed with the proof of Theorem \ref{t2}, which we split into several
steps. \\[3mm]
{\bf Step 1: Localization.} For every operator $A \in L(l^p)$, we denote
its coset $A + K(l^p)$ in the Calkin algebra by $A^\pi$, and for every
multiplier $a \in M^p$, we put $\tilde{a}(t) := a(1/t)$. The identities
\begin{equation} \label{e7}
T(ab) = T(a)T(b) + H(a)H(\tilde{b}) \quad \mbox{and} \quad H(ab) = T(a) H(b) +
H(a) T(\tilde{b}),
\end{equation}
which hold for arbitrary $a, \, b \in M^p$, together with the compactness of the
Hankel operators $H(c)$ for $c \in C_p$ show that the set $\cC_p$ of all cosets
$T(c)^\pi$ with $c \in C_p$ and $c = \tilde{c}$ forms a central subalgebra of
the algebra $\eT\eH(M^p)/K(l^p)$ and, in particular, of the algebra
$\eT\eH(PC_p)/K(l^p)$. One can, thus, reify Allan's local principle with
$\eT\eH(PC_p)/K(l^p)$ and $\cC_p$ in place of $\cA$ and $\cC$, respectively. It
is not hard to see that the maximal ideal space of $\cC_p$ is homeomorphic to
the arc $\sT_+$, with $t \in \sT_+$ corresponding to the maximal ideal $\{c \in
\cC_p: c(t) = 0\}$ of $\cC_p$. We let $\cJ_t$ denote the smallest closed ideal
of $\eT\eH(PC_p)/K(l^p)$ which contains the maximal ideal $t$ and write
$A_t^\pi$ for the coset $A^\pi + \cJ_t$ of $A \in \eT\eH(PC_p)$. Instead of
$T(a)^\pi_t$ and $H(b)^\pi_t$ we often write $T^\pi_t(a)$ and $H^\pi_t(b)$,
respectively, and the local quotient algebra $(\eT\eH(PC_p)/K(l^p))/\cJ_t$ is
denoted by $\eT\eH^\pi_t(PC_p)$ therefore. By Allan's local principle, we then
have
\begin{equation} \label{e8}
\sigma_{\eT\eH(PC_p)/K(l^p)} (A^\pi) = \cup_{t \in \sT_+}
\sigma_{\eT\eH^\pi_t(PC_p)} (A^\pi_t)
\end{equation}
for every $A \in \eT\eH(PC_p)$. \\[3mm]
{\bf Step 2: Local equivalence of multipliers.} Let $a, \, b \in PC_p$ and
$t \in \sT_+$. We show that if $a(t^\pm) = b(t^\pm)$ and $a(\overline{t}^\pm)
= b(\overline{t}^\pm)$, then $T^\pi_t(a) = T^\pi_t(b)$ and $H^\pi_t(a) =
H^\pi_t(b)$. This fact will be used in what follows in order to replace
multipliers by locally equivalent ones. It is clearly sufficient to prove that
if $a \in PC_p$ satisfies $a(t^\pm) = a(\overline{t}^\pm) = 0$, then $T^\pi(a),
\, H^\pi(a) \in \cJ_t$. We will give this proof for $t \in \sT_+^0$; the
proof for for $t = \pm 1$ is similar.

Given $\varepsilon > 0$, let $f \in P\sC$ such that $\|a - f\|_{M_p} <
\varepsilon$. Then there is an open arc $U := (e^{-i \delta} t, \, e^{i \delta}
t) \subset \sT_+$ such that $|a(s)| < \varepsilon$ almost everywhere on $U \cup
\overline{U}$ and such that $f$ has at most one discontinuity in each of $U$
and $\overline{U}$. Then $|f(s)| < 2 \varepsilon$ for $s \in U \cup
\overline{U}$. Now choose a real-valued function $\varphi_0 \in C^\infty(\sT)$
such that $\varphi_0(t) = 1$, the support of $\varphi_0$ is contained in $U$,
and $\varphi_0$ is monotonously increasing on the arc $(e^{-i \delta} t, \, t)$
and monotonously decreasing on $(t, \, e^{i \delta} t)$. Set $\varphi :=
\varphi_0 + \widetilde{\varphi_0}$. Then $\varphi = \widetilde{\varphi}$, and
\[
T^\pi(f) - T^\pi(f\varphi) = T^\pi(f(1-\varphi)) = T^\pi(f) T^\pi(1-\varphi) \in
\cJ_t,
\]
\[
H^\pi(f) - H^\pi(f\varphi) = H^\pi(f(1-\varphi)) = H^\pi(f) T^\pi(1-\varphi) \in
\cJ_t.
\]
Since $\|f\varphi\|_\infty < 2 \varepsilon$ and $\mbox{Var} (f \varphi) < 8
\varepsilon$, we conclude that $\|f \varphi\|_{M_p} < 10 c_p \varepsilon$ from
Stechkin's inequality. Thus, $\|T^\pi(f \varphi)\| < 10 c_p \varepsilon$ and
$\|H^\pi(f \varphi)\| < 10 c_p \varepsilon$, with a constant $c_p$ depending on
$p$ only. Thus, $T^\pi(a)$ differs from the element $T^\pi(f) - T^\pi(f\varphi)
\in \cJ_t$ by the element $T^\pi(a-f) + T^\pi(f\varphi)$, which has a norm less
than $(1 + 10c_p)\varepsilon$. Since $\varepsilon > 0$ is arbitrary and $\cJ_t$
is closed, this implies $T^\pi(a) \in \cJ_t$. Analogously, $H^\pi(a) \in
\cJ_t$. \\[3mm]
{\bf Step 3: The local algebras at $t \in \sT_+^0$.} We start with
describing the local algebras $\eT\eH^\pi_t(PC_p)$ at points $t \in \sT_+^0$.
Let $\chi_t$ denote the characteristic function of the arc in
$\sT$ which connects $t$ with $\bar{t}$ and runs through the point -1. Clearly,
$\chi_t \in PC_p$. The crucial observation, which is a simple consequence of the
identities (\ref{e7}), is that the operator $T(\chi_t) + H(\chi_t)$ is an
idempotent. Further, let $\varphi_t \in C_p$ be any multiplier such that $0 \le
\varphi_t \le 1$, $\varphi_t(t) = 1$, $\varphi(\bar{t}) = 0$ and $\varphi_t +
\widetilde{\varphi_t} = 1$. Again by (\ref{e7}), the coset $T^\pi_t(\varphi_t)$
is an idempotent.

We claim that the idempotents $p_t := T^\pi_t(\varphi_t)$ and $q_t := T^\pi_t
(\chi_t) + H^\pi_t (\chi_t)$ together with the identity element $e := I^\pi_t$
generate the local algebra $\eT\eH^\pi_t(PC_p)$. Let $a, \, b \in PC_p$. Then,
using step 2,
\begin{eqnarray} \label{e9}
T^\pi_t(a) & = &  a(t^+) T^\pi_t (\chi_t \varphi_t) + a(t^-) T^\pi_t ((1-\chi_t)
\varphi_t) + a(\bar{t}^-) T^\pi_t (\chi_t (1 - \varphi_t)) \nonumber\\
& & \quad + \; a(\bar{t}^+) T^\pi_t ((1-\chi_t)(1-\varphi_t)).
\end{eqnarray}
It is not hard to check that
\begin{eqnarray} \label{e10}
T^\pi_t (\chi_t \varphi_t) & = & p_t q_t p_t, \nonumber \\
T^\pi_t ((1-\chi_t) \varphi_t) & = & p_t (e-q_t) p_t, \nonumber \\
T^\pi_t (\chi_t (1 - \varphi_t)) & = & (e-p_t) q_t (e-p_t), \nonumber \\
T^\pi_t ((1-\chi_t)(1-\varphi_t)) & = & (e-p_t) (e-q_t) (e-p_t).
\end{eqnarray}
Let us verify the first of these identities, for example. By definition,
\[
p_t q_t p_t = T^\pi_t(\varphi_t) T^\pi_t (\chi_t) T^\pi_t(\varphi_t) +
T^\pi_t(\varphi_t) H^\pi_t (\chi_t) T^\pi_t(\varphi_t).
\]
Since $T(\varphi_t)$ commutes with $T(\chi_t)$ modulo compact operators and
$H(\widetilde{\varphi_t})$ is compact, we can use the identities (\ref{e7}) to
conclude
\[
T^\pi_t(\varphi_t) T^\pi_t (\chi_t) T^\pi_t(\varphi_t) = T^\pi_t (\chi_t)
T^\pi_t(\varphi_t) = T^\pi_t (\chi_t \varphi_t).
\]
Further, due to the compactness of $H(\varphi_t)$ and
$H(\widetilde{\varphi_t})$,
\[
T^\pi_t(\varphi_t) H^\pi_t (\chi_t) T^\pi_t(\varphi_t) = H^\pi_t (\varphi_t
\chi_t) T^\pi_t(\varphi_t) = H^\pi_t (\varphi_t \chi_t \widetilde{\varphi_t}).
\]
Since $\varphi_t \chi_t \widetilde{\varphi_t}$ is a continuous function,
$H^\pi_t (\varphi_t \chi_t \widetilde{\varphi_t}) = 0$. This gives the first of
the identities (\ref{e10}). The others follow in a similar way. Thus, (\ref{e9})
and (\ref{e10}) imply that $T^\pi_t(a)$ belongs to the algebra generated by $e$,
$p_t$ and $q_t$. Similarly, we write
\begin{eqnarray} \label{e11}
H^\pi_t(b) & = & b(t^+) H^\pi_t (\chi_t \varphi_t) + b(t^-) H^\pi_t ((1-\chi_t)
\varphi_t) + b(\bar{t}^-) H^\pi_t (\chi_t (1 - \varphi_t)) \nonumber \\
& & \quad + \; b(\bar{t}^+) H^\pi_t ((1-\chi_t)(1-\varphi_t))
\end{eqnarray}
and use the identities
\begin{eqnarray} \label{e12}
H^\pi_t (\chi_t \varphi_t) & = & p_t q_t (e-p_t), \nonumber \\
H^\pi_t ((1-\chi_t) \varphi_t) & = & - p_t q_t (e-p_t), \nonumber \\
H^\pi_t (\chi_t (1 - \varphi_t)) & = & (e-p_t)q_t p_t, \nonumber \\
H^\pi_t ((1-\chi_t)(1-\varphi_t)) & = & - (e-p_t) q_t p_t
\end{eqnarray}
to conclude that $H^\pi_t(b)$ also belongs to the algebra generated by $e$,
$p_t$ and $q_t$. Thus, the algebra $\eT\eH^\pi_t(PC_p)$ is subject to the two
idempotents theorem.

In order to apply this theorem we have to determine the spectrum of the coset
$p_t q_t p_t = T^\pi_t (\chi_t \varphi_t)$ in that algebra. We claim that
\begin{equation} \label{e13}
\sigma_{\eT\eH^\pi_t(PC_p)} (T^\pi_t (\chi_t \varphi_t)) = \{\mu_q(\lambda) :
\lambda \in \overline{\sR} \}
\end{equation}
with $1/p + 1/q = 1$. Let $a_t \in PC_p$ be a multiplier with the following
properties: \\[1mm]
\hspace*{3mm} $(a)$ $a_t$ is continuous on $\sT \setminus \{t\}$ and has a jump
at $t \in \sT$. \\[1mm]
\hspace*{3mm} $(b)$ $a_t(t^+) = \chi_t(t^+) = 1$ and $a_t(t^-) = \chi_t(t^-) =
0$. \\[1mm]
\hspace*{3mm} $(c)$ $a_t$ takes values in $\{\mu_q(\lambda) : \lambda \in
\overline{\sR} \}$ only. \\[1mm]
\hspace*{3mm} $(d)$ $a_t$ is zero on the arc joining $-t$ to $t$ which contains
the point 1. \\[1mm]
Then, by Theorem \ref{t1}, the essential spectrum of the Toeplitz operator
$T(a_t)$ in each of the algebras $L(l^p)/K(l^p)$ and $\eT(PC_p)/K(l^p)$ is
equal to the arc $\{\mu_q(\lambda) : \lambda \in \overline{\sR} \}$. Hence, the
essential spectrum of $T(a_t)$, now considered as an element of the algebra
$\eT\eH(PC_p)/K(l^p)$, is also equal to this arc. Hence,
\[
\sigma_{\eT\eH^\pi_t(PC_p)} (T^\pi_t (a_t)) \subseteq \{\mu_q(\lambda) :
\lambda \in \overline{\sR} \}
\]
by Allan's local principle. Since $T^\pi_t(a_t) = T^\pi_t (\chi_t \varphi_t)$,
this settles the inclusion $\subseteq$ in (\ref{e13}). For the reverse
inclusion, let $b_t\in PC_p$ be a multiplier with the following
properties: \\[1mm]
\hspace*{3mm} $(a)$ $b_t$ is continuous on $\sT \setminus \{t\}$ and has a jump
at $t \in \sT$. \\[1mm]
\hspace*{3mm} $(b)$ $b_t(t^\pm) = \chi_t(t^\pm)$. \\[1mm]
\hspace*{3mm} $(c)$ $b_t$ takes values not in $\{\mu_q(\lambda) : \lambda \in
\sR \}$ on the arc joining $-t$ to $t$ which \\
\hspace*{3mm} \hphantom{$(c)$} contains the point $-1$. \\[1mm]
\hspace*{3mm} $(d)$ $b_t$ is zero on the arc joining $-t$ to $t$ which contains
the point 1. \\[1mm]
Then, again by Theorem \ref{t1}, the essential spectrum of the Toeplitz operator
$T(b_t)$ in each of the algebras $L(l^p)/K(l^p)$ and $\eT(PC_p)/K(l^p)$ is
equal to the union of the arc $\{\mu_q(\lambda) : \lambda \in \overline{\sR}
\}$ and the range of $b_t$. Hence, the essential spectrum of $T(b_t)$, now
considered as an element of the algebra $\eT\eH(PC_p)/K(l^p)$, is also equal to
this union. Since $b_t$ is continuous on $\sT \setminus \{t\}$ by property
$(a)$, we have
\[
\sigma_{\eT\eH^\pi_s(PC_p)} (T^\pi_s (b_t)) = \{ b_t(s), \, b_t(\bar{s}) \}
\]
for $s \in \sT^0_+ \setminus \{t\}$. Since the points $b_t(s)$ and
$b_t(\bar{s})$ do not belong to $\{\mu_q(\lambda) : \lambda \in
\sR \}$ by property $(c)$, we conclude that the open arc $\{\mu_q(\lambda) :
\lambda \in \sR \}$ is contained in the local spectrum of $T(b_t)$ at $t$.
Since spectra are closed, this implies
\[
\{\mu_q(\lambda) : \lambda \in \overline{\sR} \} \subseteq
\sigma_{\eT\eH^\pi_t(PC_p)} (T^\pi_t (b_t)).
\]
Since $T^\pi_t(b_t) = T^\pi_t (\chi_t \varphi_t)$ by property $(b)$,
this settles the inclusion $\supseteq$ in (\ref{e13}).

Since $\nu_q(\lambda)^2 = \mu_q(\lambda) (1-\mu_q(\lambda))$, we can choose
$\sqrt{\mu_q(\lambda) (1-\mu_q(\lambda))} = \nu_q(\lambda)$. With this choice
and identities (\ref{e9}) -- (\ref{e12}) it becomes evident that the two
idempotents theorem associates with the coset $T^\pi_t(a) + H^\pi_t(b)$ the
matrix function
\[
\lambda \mapsto \pmatrix{a(t^+) \mu_q(\lambda) + a(t^-)(1-\mu_q(\lambda)) &
(b(t^+)-b(t^-)) \nu_q(\lambda) \cr
(b(\bar{t}^-)-b(\bar{t}^+)) \nu_q(\lambda) &
a(\bar{t}^-) (1-\mu_q(\lambda)) + a(\bar{t}^+) \mu_q(\lambda)}
\]
on $\overline{\sR}$. \\[3mm]
{\bf Step 4: The local algebra at $1 \in \sT_+$.} Next we are going to consider
the local algebra $\eT\eH^\pi_1(PC_p)$ at the fixed point 1 of the mapping $t
\mapsto \bar{t}$. Let $f : \sT \to \sC$ denote the function $e^{is} \mapsto 1
-s/\pi$ where $s \in [0, 2 \pi)$. This function belongs to $PC_p$, and it has
its only jump at the point $1 \in \sT$ where $f(1^\pm) = \pm 1$. Using ideas
from \cite{PR1}, it was shown in \cite{Ro1} by one of the authors that the
Hankel operator $H(f)$ belongs to the Toeplitz algebra $\eT(PC_p)$ and that its
essential spectrum is given by
\begin{equation} \label{e14}
\sigma_{ess} (H(f)) = \{ 2i \, \nu_q(\lambda) : \lambda \in \overline{\sR} \}.
\end{equation}
(in fact, this identity was derived in \cite{Ro1} with $p$ in place of $q$,
which makes no difference since $\nu_p(- \lambda) = \nu_q(\lambda)$ for
every $\lambda$.) Let $\chi_+$ denote the characteristic function of the upper
half-circle $\sT_+$. Since every coset $T^\pi_1(a)$ with $a \in PC_p$ is a
linear combination of the cosets $I^\pi_1$ and $T^\pi_1(\chi_+)$ and every coset
$H^\pi_1(b)$ is a multiple of the coset $H^\pi_1(f)$, the local algebra
$\eT\eH^\pi_1(PC_p)$ is singly generated (as a unital algebra) by the coset
$T^\pi_1(\chi_+)$. In particular, $\eT\eH^\pi_1(PC_p)$ is a commutative Banach
algebra, and its maximal ideal space is homeomorphic to the spectrum of its
generating element. Similar to the proof of (\ref{e13}) one can show that
\begin{equation} \label{e15}
\sigma_{\eT\eH^\pi_1(PC_p)} (T^\pi_t (\chi_+)) = \{\mu_q(\lambda) :
\lambda \in \overline{\sR} \}
\end{equation}
It is convenient for our purposes to identify the maximal ideal space of the
algebra $\eT\eH^\pi_1(PC_p)$ with $\overline{\sR}$. The Gelfand transform of
$T^\pi_t (\chi_+)$ is then given by $\lambda \mapsto \mu_q(\lambda)$ due to
(\ref{e15}). Let $h$ denote the Gelfand transform of $H^\pi_1(f)$. From
(\ref{e7}) we obtain
\[
H^\pi_1(f)^2 = T^\pi_1(f \tilde{f}) - T^\pi_1(f) T^\pi_1(\tilde{f}).
\]
The function $f \tilde{f}$ is continuous at $1 \in \sT$ and has the value
$-1$ there, and the function $f + \tilde{f}$ is continuous at $1 \in \sT$ and
has the value $0$ there. Thus,
\[
H^\pi_1(f)^2 = - I^\pi_1 + T^\pi_1(f)^2.
\]
Since $T^\pi_1(f) = T^\pi_1(2 \chi_+ -1) = 2 T^\pi_1(\chi_+) - I^\pi_1$ we
conclude that
\[
h(\lambda)^2 = (2 \mu_q(\lambda) -1)^2 -1 = (\sinh (\pi (\lambda + i/q)))^{-2}
\]
if $\lambda \in \sR$ and by $h(\pm \infty) = 0$. By (\ref{e14}), this equality
necessarily implies that
\[
h(\lambda) = (\sinh (\pi (\lambda + i/q)))^{-1} = 2i \nu_q(\lambda)
\]
if $\lambda \in \sR$ and $h(\pm \infty) = 0$. Combining these results we
find that the Gelfand transform of $T^\pi_1(a) + H^\pi_1(b)$ is the function
\[
\lambda \mapsto a(1^+) \mu_q(\lambda) + a(1^-) (1 -  \mu_q(\lambda)) + i \,
(b(1^+) - b(1^-)) \nu_q (\lambda).
\]
{\bf Step 5: The local algebra at $-1 \in \sT_+$.}
It remains to examine the local algebra $\eT\eH^\pi_{-1}(PC_p)$ at the point
$-1$. Let $\Lambda : l^2 \to l^2$ denote the mapping $(x_n)_{n \ge 0} \mapsto
((-1)^n x_n)_{n \ge 0}$. Clearly, $\Lambda^{-1} = \Lambda$, and one easily
checks (perhaps most easily on the level of the matrix entries, which are
Fourier coefficients) that
\[
\Lambda^{-1} T(a) \Lambda = T(\hat{a}) \quad \mbox{and} \quad
\Lambda^{-1} H(a) \Lambda = - H(\hat{a})
\]
for $a \in PC_p$, where $\hat{a}(t) := a(-t)$. Thus, the mapping $A \mapsto
\Lambda^{-1} A \Lambda$ is an automorphism of the algebra $\eT\eH(PC_p)$, which
maps compact operators to compact operators and induces, thus, an automorphism
of the algebra $\eT\eH(PC_p)/K(l^p)$. The latter maps the local ideal at $1$ to
the local ideal at $-1$ and vice versa and induces, thus, an isomorphism
between the local algebras $\eT\eH^\pi_1(PC_p)$ and $\eT\eH^\pi_{-1}(PC_p)$,
which sends $T^\pi_1(\chi_+)$ to $T^\pi_{-1} (1-\chi_+)$ and $H^\pi_1(\chi_+)$
to $- H^\pi_{-1} (1-\chi_+) = H^\pi_{-1} (\chi_+)$, respectively. \\[3mm]
{\bf Step 6: From local to global invertibility.}
We have identified the right-hand sides of (\ref{e3}) and (\ref{e4}) as the
functions which are locally associated with the operator $T(a) + H(b)$ via the
two idempotents theorem and via Gelfand theory for commutative Banach algebras,
respectively. It follows from the two idempotents theorem and from Gelfand
theory that the so-defined mappings $\smb (t, \, \lambda)$ extend to a
continuous homomorphism from $\eT\eH(PC_p)$ to $\sC^{2 \times 2}$ or $\sC$,
respectively, which combine to a continuous homomorphism from $\eT\eH(PC_p)$
to the algebra $\cF$. Allan's local principle then implies that the coset $A +
K(l^p)$ of an operator $A \in \eT\eH(PC_p)$ is invertible in
$\eT\eH(PC_p)/K(l^p)$ if and only if its symbol does not vanish. The proof of
estimate (\ref{e241111.1}) will base on Mellin homogenization arguments. We
therefore postpone it until Section \ref{s5}; see estimate (\ref{e27}).
\\[3mm]
{\bf Step 7: Inverse closedness.}
It remains to show that $\eT\eH(PC_p)/K(l^p)$ is an inverse closed subalgebra of
the Calkin algebra $L(l^p)/K(l^p)$. We shall prove this fact by using a {\em
thin spectra argument} as follows: If $\cA$ is a unital closed subalgebra of a
unital Banach algebra $\cB$, and if the spectrum in $\cA$ of every element in a
dense subset of $\cA$ is thin, i.e. if its interior with respect to the topology
of $\sC$ is empty, then $\cA$ is inverse closed in $\cB$. See, e.g.,
\cite{RSS1}, Corollary 1.2.32, for a simple proof of this argument.

Let $\cA_0$ be the set of all operators of the form
\begin{equation} \label{e16}
A := \sum_{i=1}^l \prod_{j=1}^k (T(a_{ij}) + H(b_{ij})) \quad \mbox{with} \;
A_{ij}, \, b_{ij} \in P\sC,
\end{equation}
and write $\sigma^{TH}_{ess} (A)$ for the spectrum of $A$ in
$\eT\eH(PC_p)/K(l^p)$. Then $\cA_0/K(l^p)$ is dense in $\eT\eH(PC_p)/K(l^p)$,
and the assertion will follow once we have shown that $\eT\eH(PC_p)/K(l^p)$ is
thin for every $A \in \cA_0$.

Given $A$ of the form (\ref{e16}), let $\Omega$ denote the set of all
discontinuities of the functions $a_{ij}$ and $b_{ij}$, and put
$\widetilde{\Omega} := (\Omega \cup \overline{\Omega}) \cap \sT_+$. Clearly,
$\widetilde{\Omega}$ is a finite set. By what we have shown above,
\[
\sigma^{TH}_{ess} (A) = \cup_{(t,\lambda) \in \sT_+ \times \overline{\sR}} \,
\sigma (\smb (A) (t, \, \lambda))
\]
where $\sigma (B)$ stands for the spectrum (= set of the eigenvalues) of the
matrix $B$. We write $\sigma^{TH}_{ess} (A)$ as $\Sigma_1 \cup \Sigma_2 \cup
\Sigma_3$ where
\begin{eqnarray*}
\Sigma_1 & := & \cup_{(t,\lambda) \in \{-1,1\} \times \overline{\sR}} \,
\sigma (\smb (A) (t, \, \lambda)), \\
\Sigma_2 & := & \cup_{(t,\lambda) \in (\sT^0_+ \setminus \widetilde{\Omega})
\times \overline{\sR}} \, \sigma (\smb (A) (t, \, \lambda)), \\
\Sigma_3 & := & \cup_{(t,\lambda) \in (\widetilde{\Omega} \setminus \{-1,1\})
\times \overline{\sR}} \, \sigma (\smb (A) (t, \, \lambda)).
\end{eqnarray*}
It is clear that $\Sigma_1$ is a set of measure zero. It is also clear that each
set
\[
\Sigma_{2,t} := \cup_{\lambda \in \overline{\sR}} \, \sigma (\smb (A) (t,
\, \lambda)) \quad \mbox{with} \; t \in \sT^0_+ \setminus \widetilde{\Omega}
\]
has measure zero. Since the functions $a_{ij}$ and $b_{ij}$ are piecewise
constant, the mapping $t \mapsto \Sigma_{2,t}$ is constant on each connected
component of $\sT^0_+ \setminus \widetilde{\Omega}$, and the number of
components is finite. Thus, $\Sigma_2$ is actually a finite union of sets of
measure zero. Since $\widetilde{\Omega}$ is finite, it remains to show that each
of the sets
\[
\Sigma_{3,t} := \cup_{\lambda \in \overline{\sR}} \, \sigma (\smb (A) (t, \,
\lambda))  \quad \mbox{with} \; t \in \widetilde{\Omega} \setminus \{-1,1\}
\]
has measure zero. For this goal it is clearly sufficient to show that each set
\[
\Sigma_{3,t}^0 := \cup_{\lambda \in \sR} \, \sigma (\smb (A) (t, \,
\lambda))  \quad \mbox{with} \; t \in \widetilde{\Omega} \setminus \{-1,1\}
\]
has measure zero. Let $t \in \widetilde{\Omega} \setminus \{-1,1\}$, and write
$\smb (A) (t, \, \lambda)$ as $(c_{ij}(\lambda))_{i,j = 1}^2$. The eigenvalues
of this matrix are $s_\pm (\lambda) = (c_{11}(\lambda) + c_{22}(\lambda))/2 \pm
\sqrt{r(\lambda)}$ where
\[
r(\lambda) = (a_{11}(\lambda) + a_{22}(\lambda))^2/4 - (a_{11}(\lambda)
a_{22}(\lambda) - a_{12}(\lambda) a_{21}(\lambda))
\]
and where $\sqrt{r(\lambda)}$ is any complex number the square of which is
$r(\lambda)$. Since $r$ is composed by the meromorphic functions $\coth$ and
$1/\sinh$, the set of zeros of $r$ is discrete. Hence, $\sR \setminus \{
\lambda \in \sR : r(\lambda) = 0 \}$ is an open set, which as the union of an
at most countable family of open intervals. Let $I$ be one of these intervals.
Then $I$ can be represented as the union of countably many compact subintervals
$I_n$ such that the intersection $I_n \cap I_m$ consists of at most one
point whenever $n \neq m$ and each set $r(I_n)$ is contained in a domain where
a continuous branch, say $f_n$, of the function $z \mapsto \sqrt{z}$ exists.
Then $\pm f_n \circ r : I_n \to \sC$ is a continuously differentiable
function, which implies that $(\pm f_n \circ r)(I_n)$ is a set of measure zero.
Consequently, the associated sets $s_\pm(I_n)$ of eigenvalues have measure
zero, too. Since the countable union of sets of measure zero has measure zero,
we conclude that each set $\Sigma_{3,t}^0$ has measure zero, which finally
implies that $\sigma^{TH}_{ess} (A) = \Sigma_1 \cup \Sigma_2 \cup
\Sigma_3$ has measure zero and is, thus, thin. This settles the proof of
the inverse closedness and concludes the proof of Theorem \ref{t2}. \hfill \qed
\\[3mm]
We would like to mention that there is another proof of the inverse closedness
assertion in the previous theorem which is based on ideas from \cite{FRS} and
which works also in other situations.  
\section{An extended Toeplitz algebra} \label{s3}
In the proof of the announced index formula for Toeplitz plus Hankel operators,
we shall need an extension of the results of the previous section to certain
matrix operators. For $k \in \sN$ and $X$ a linear space, we let $X_k$ and
$X_{k \times k}$ stand for the linear spaces of all vectors of length $k$ and
of all $k \times k$-matrices with entries in $X$, respectively. If $X$ is an
algebra, then $X_{k \times k}$ becomes an algebra under the standard matrix
operations. If $X$ is a Banach space, then $X_k$ and $X_{k \times k}$ become
Banach spaces with respect to the norms
\begin{equation} \label{ne1}
\| (x_j)_{j=1}^k \| = \sum_{j=1}^k \|x_j\| \quad \mbox{and} \quad
\| (a_{ij})_{i,j = 1}^k \| = k \sup_{1 \le i,j \le k} \|a_{ij}\|.
\end{equation}
If, moreover, $X$ is a Banach algebra, then $X_{k \times k}$ is a Banach
algebra with respect to the introduced norm. Actually, any other norm on $X_k$
and any other compatible matrix norm on $X_{k \times k}$ will do the same job.
Note also that if $X$ is a $C^*$-algebra there is a unique norm (different from
the above mentioned) which makes $X_{k \times k}$ to a $C^*$-algebra. Since we
will not employ $C^*$-arguments, the choice (\ref{ne1}) will be sufficient for
our purposes.

Let $\eT^0(PC_p)$ denote the smallest closed subalgebra of $L(l^p(\sZ))$ which
contains the projection $P$ and all Laurent operators $L(a)$ with $a \in PC_p$.
The algebra $\eT^0(PC_p)$ contains $\eT(PC_p)$ in the sense that the operator
$P L(a) P : \im P \to \im P$ can be identified with the Toeplitz operator
$T(a)$. For $k \in \sN$, the matrix algebra $\eT^0(PC_p)_{k \times k}$
will be also denoted by $\eT^0_{k \times k}(PC_p)$. One can characterize
$\eT^0_{k \times k}(PC_p)$ also as the smallest closed subalgebra of
$L(l^p(\sZ)_k)$ which contains all operators of the form $L(a) \diag P + L(b)
\diag Q$ with $a, \, b \in (PC_p)_{k \times k}$, where $Q := I - P$, $\diag A$
stands for the operator on $L(l^p(\sZ)_k)$ which has $A \in L(l^p(\sZ))$ at
each entry of its main diagonal and zeros at all other entries, and where $L(a)
= (L(a_{ij}))_{i,j = 1}^k$ refers to the matrix Laurent operator with generating
function $a = (a_{ij})_{i,j = 1}^k$. Note that $K(l^p(\sZ)_k)$ is contained in
$\eT^0_{k \times k}(PC_p)$.

The Fredholm theory for operators in $\eT^0_{k \times k}(PC_p)$ is well known.
We will present it in a form which is convenient for our purposes. Our main
tools are again Allan's local principle (Theorem \ref{t5}) and a matrix version
of the two idempotents theorem (Theorem \ref{t6}) due to \cite{FRS}. Here is
the result.
\begin{theo} \label{nt2}
Let $a, \, b \in (PC_p)_{k \times k}$. \\[1mm]
$(a)$ The operator $A := L(a) \diag P + L(b) \diag Q$ is Fredholm on
$l^p(\sZ)_k$ if and only if the matrix
\begin{eqnarray*}
\lefteqn{(\smb A)(t, \, \lambda) =} \\
&& \hspace*{-5mm} \pmatrix{
a(t^-) + (a(t^+) - a(t^-)) \diag \mu_q(\lambda) &
(b(t^+) - b(t^-)) \diag \nu_q(\lambda) \cr
(a(t^+) - a(t^-)) \diag \nu_q(\lambda) &
b(t^+) - (b(t^+) - b(t^-)) \diag \mu_q(\lambda)}
\end{eqnarray*}
is invertible for every pair $(t, \, \lambda) \in \sT \times \overline{\sR}$.
\\[1mm]
$(b)$ The mapping $\smb$ defined in assertion $(a)$ extends to a continuous
algebra homomorphism from $\eT^0_{k \times k}(PC_p)$ to the algebra $\cF$ of
all bounded functions on $\sT \times \overline{\sR}$ with values in $\sC_{2k
\times 2k}$. Moreover, there is a constant $M$ such that
\begin{equation} \label{e271011.2}
\| \smb A \| := \sup_{(t, \lambda) \in \sT_+ \times \overline{\sR}} \|\smb A 
(t, \, \lambda)\|_\infty \le M \inf_{K \in K(l^p(\sZ)_k)} \|A + K\|
\end{equation}
for every operator $A \in \eT^0_{k \times k}(PC_p)$. \\[1mm]
$(c)$ An operator $A \in \eT^0_{k \times k}(PC_p)$ has the Fredholm property on
$l^p(\sZ)_k$ if and only if the function $\smb A$ is invertible in $\cF$.
\\[1mm]
$(d)$ The algebra $\eT^0_{k \times k}(PC_p)/K(l^p(\sZ)_k)$ is inverse closed in
the Calkin algebra $L(l^p(\sZ)_k)/K(l^p(\sZ)_k)$. \\[1mm]
$(e)$ If $A \in \eT^0_{k \times k}(PC_p)$ is a Fredholm operator, then
\[
\ind A = - \wind (\det \smb A (t, \, \lambda) /(\det a_{22}(t, \, \infty)
\det a_{22}(t, \, -\infty)))
\]
where $\smb A = (a_{ij})_{i,j = 1}^2$ with $k \times k$-matrix-valued functions
$a_{ij}$.
\end{theo}
It is a non-trivial fact that the function
\[
W : \sT \times \overline{\sR}, \quad (t, \, \lambda) \mapsto \det \smb A (t, \,
\lambda) /(\det a_{22}(t, \, \infty) \det a_{22}(t, \, -\infty))
\]
forms a closed curve in the complex plane. Thus, the winding number of $W$ is
well defined if $A$ is a Fredholm operator.

The remainder of this section is devoted to the proof of Theorem \ref{nt2}. We
shall mainly make use of results from Sections 2.3 - 2.5 in \cite{HRS1} and
Chapter 6 in \cite{BS1}. We will be quite sketchy when the arguments are
close to those from the proof of Theorem \ref{t2}. \\[3mm]
{\bf Step 1: Spline spaces.} We start with recalling some facts about spline
spaces and operators thereon from \cite{HRS1}. Let $\chi_{[0,1]}$ denote the
characteristic function of the interval $[0, \, 1] \subset \sR$ and, for $n \in
\sN$, let $S_n$ denote the smallest closed subspace of $L^p(\sR)$ which contains
all functions
\[
\varphi_{k,n}(t) := \chi_{[0,1]}(nt-k), \quad t \in \sR,  
\]
where $k \in \sZ$. The space $l^p(\sZ)$ can be identified with each of the
spaces $S_n$ in the sense that a sequence $(x_k)$ is in $l^p(\sZ)$ if and only
if the series $\sum_{k \in \sZ} x_k \varphi_{k,n}$ converges in $L^p(\sR)$ and
that
\[
\left\| \sum x_k \varphi_{k,n} \right\|_{L^p(\sR)} = n^{-1/p}
\left\| (x_k) \right\|_{l^p(\sZ)}
\]
in this case. Thus, the linear operator
\[
E_n : l^p(\sZ) \to S_n \subset L^p(\sR), \quad (x_k) \mapsto n^{1/p} \sum x_k
\varphi_{k,n}, 
\]
and its inverse $E_{-n} : L^p(\sR) \supset S_n \to l^p(\sZ)$ are isometries
for every $n$. Further we define operators 
\[
L_n : L^p(\sR) \to S_n, \, \quad u \mapsto n \sum_{k \in \sZ} \langle u, \,
\varphi_{k,n} \rangle \varphi_{k,n}
\]
with respect to the sesqui-linear form $\langle u, \, v \rangle := \int_\sR u
\overline{v} dx$, where $u \in L^p(\sR)$ and $v \in L^q(\sR)$ with $1/p + 1/q =
1$. It is easy to see that every $L_n$ is a projection operator with norm 1 and
that the $L_n$ converge strongly to the identity operator on $L^p(\sR)$ as $n
\to \infty$. Finally we set  
\[
Y_t : l^p(\sZ) \to l^p(\sZ), \quad (x_k) \mapsto (t^{-k} x_k) \quad
\mbox{for} \; t \in \sT.
\]
Clearly, $Y_t$ is an isometry, and $Y_t^{-1} = Y_{t^{-1}}$. One easily checks
that $Y_t^{-1} L(a) Y_t = L(a_t)$ with $a_t(s) = a(ts)$ for every multiplier
$a$, which implies in particular that $Y_t^{-1} \eT^0(PC_p) Y_t = \eT^0(PC_p)$.
\\[3mm]
{\bf Step 2: Some homomorphisms.} In Sections 2.3.3 and 2.5.2 of \cite{HRS1} it
is shown that, for every $A \in \eT^0(PC_p)$ and every $t \in \sT$, the strong
limit
\[
\mbox{smb}_t A := \mbox{s-lim}_{n \to \infty} E_n Y_t^{-1} A Y_t E_{-n} L_n
\]
exists and that the mapping $\mbox{\rm smb}_t$ is a bounded unital algebra
homomorphism. This homomorphism can be extended in a natural way to the matrix
algebra $\eT^0_{k \times k}(PC_p)$. We denote this extension by
$\mbox{smb}_t A$ again.

In order to characterize the range of the homomorphism $\mbox{\rm smb}_t$, we
have to introduce some operators on $L^p(\sR)$. Let $\chi_+$ stand for the
characteristic function of the interval $\sR^+ = [0, \, \infty)$ and $\chi_+ I$
for the operator of
multiplication by $\chi_+$. Further, $S_\sR$ refers to the singular integral
operator
\[
(S_\sR f)(t) := \frac{1}{\pi i} \int_{-\infty}^\infty \frac{f(s)}{s-t} \, ds,
\]
with the integral understood as a Cauchy principal value. Both $\chi_+I$ and
$S_\sR$ are bounded on $L^p(\sR)$, and $S_\sR^2 = I$. Thus, the operators  
$P_\sR := (I+S_\sR)/2$ and $Q_\sR := I - P_\sR$ are bounded projections on
$L^p(\sR)$. We let $\Sigma_k^p (\sR)$ stand for the smallest closed
subalgebra of $L(L^p(\sR)_k)$ which contains the operators $\diag \chi_+I$,
$\diag  S_\sR$, and all operators of multiplication by constant $k \times
k$-matrix-valued functions. 
\begin{theo} \label{t271011.1}
Let $t \in \sT$. Then \\ 
$(a)$ $\mbox{\rm smb}_t \, \diag P = \diag \chi_+ I$. \\
$(b)$ $\mbox{\rm smb}_t \, L(a) = a(t^+) \diag Q_\sR + a(t^-) \diag P_\sR$ for
$a \in (PC_p)_{k \times k}$. \\                                      
$(c)$ $\mbox{\rm smb}_t \, K = 0$ for every compact operator $K$. \\
$(d)$ $\mbox{\rm smb}_t$ maps the algebra $\eT^0_{k \times k}(PC_p)$ onto
$\Sigma_k^p (\sR)$. \\
$(e)$ The algebra $\Sigma_k^p (\sR)$ is inverse closed in $L(L^p(\sR)_k)$.
\end{theo}
Assertion $(c)$ of the previous theorem implies that every mapping $\mbox{\rm
smb}_t$ induces a natural quotient homomorphism from $\eT^0(PC_p)/K(l^p(\sZ))$
to $\Sigma_1^p(\sR)$. We denote this quotient homomorphism by $\mbox{smb}_t$
again. It now easily seen that the estimate (\ref{e271011.2}) holds for every $A
\in \eT^0_{k \times k}(PC_p)$ (with the constant $M=1$ for $k=1$). \\[3mm]
{\bf Step 3: The Fredholm property.} Since the commutator $L(a) P -
P L(a)$ is compact for every $a \in C_p$, the algebra $\cC_p := \{\diag L(a) : a
\in C_p\}/K(l^p(\sZ)_k)$ lies in the center of the algebra $\cA := \eT^0_{k
\times k}(PC_p)/K(l^p(\sZ)_k)$. It is not hard to see that $\cC_p$ is isomorphic
to $C_p$; hence the maximal ideal space of $\cC_p$ is homeomorphic to the unit
circle $\sT$. In accordance with Allan's local principle, we introduce the local
ideals $\cJ_t$ and the local algebras $\cA_t := \cA/\cJ_t$ at $t \in \sT$.

By Theorem \ref{t271011.1} $(b)$, the local ideal $\cJ_t$ lies in the kernel of
$\mbox{smb}_t$. We denote the related quotient homomorphism by $\mbox{smb}_t$
again. Thus, $\mbox{smb}_t$ is an algebra homomorphism from $\cA_t$ onto
$\Sigma_k^p (\sR)$, which sends the local cosets containing the operators $\diag
P$ and $L(a)$ with $a \in (PC_p)_{k \times k}$ to $\diag \chi_+ I$ and $a(t^+)
\, \diag Q_\sR + a(t^-) \, \diag P_\sR$, respectively. By Theorem 2.3 in
\cite{HRS1}, this homomorphism is injective, i.e., it is an isomorphism between
$\cA_t$ and $\Sigma_k^p (\sR)$. 

Since $P_\sR$ and $\diag \chi_+ I$ are projections, the algebra $\Sigma_k^p
(\sR)$ is subject to the two projections theorem with coefficients, as derived
in \cite{FRS}. Alternatively, this algebra can be described by means of the
Mellin symbol calculus, see Section 2.1 in \cite{HRS1}. In each case, the result
is that an operator of the form
\begin{equation} \label{ne3}
(a^+ \diag \chi_+ I + a^- \diag \chi_- I) \, \diag P_\sR + 
(b^+ \diag \chi_+ I + b^- \diag \chi_- I)  \, \diag Q_\sR
\end{equation}
where $\chi_- := 1 - \chi_+$ and $a^\pm, \, b^\pm \in \sC_{k \times k}$ is
invertible if and only if the $(2k) \times (2k)$-matrix-valued function
\[
\lambda \mapsto \pmatrix{
a^+ \diag (1-\mu_p(\lambda)) + a^- \diag \mu_p(\lambda) & (b^+ - b^-) \, \diag
\nu_p(\lambda) \cr
(a^+ - a^-) \, \diag \nu_p(\lambda) & b^+ \diag \mu_p(\lambda) + b^- \diag
(1-\mu_p(\lambda))}
\]
is invertible at each point $\lambda \in \overline{\sR}$. Note that the function
\[
\lambda \mapsto a^+ \diag (1-\mu_p(\lambda)) + a^- \diag \mu_p(\lambda)
\] 
is continuous on $\overline{\sR}$ and that this function connects $a^+$ with
$a^-$ if $\lambda$ runs from $- \infty$ to $+ \infty$. For the sake of index
computation, one would prefer to work with a function which connects $a^-$
with $a^+$ if $\lambda$ increases. Since $\mu_p(-\lambda) = 1 - \mu_q(\lambda)$
and $\nu_p(-\lambda) = \nu_q(\lambda)$ with $q$ satisfying $1/p + 1/q = 1$, we
obtain that the operator $A$ in (\ref{ne3}) is invertible if and only if the
matrix function
\[
\lambda \mapsto \pmatrix{
a^+ \diag \mu_q(\lambda) + a^- \diag (1 - \mu_q(\lambda)) & (b^+ - b^-) \,
\diag \nu_q(\lambda) \cr
(a^+ - a^-) \, \diag \nu_q(\lambda) & b^+ \diag (1 - \mu_q(\lambda)) + b^- \diag
\mu_q(\lambda)}
\]
is invertible on $\overline{\sR}$. This observation, together with the local
principle, implies that the coset $L(a) \diag P + L(b) \diag Q +
K(l^p(\sZ)_k)$ is invertible in the quotient algebra $\eT^0_{k \times
k}(PC_p)/K(l^p(\sZ)_k)$ if and only if the matrix function in assertion $(a)$
of Theorem \ref{nt2} is invertible. In particular, this gives the ``if''-part
of assertion $(a)$. The ``only if''-part of this assertion follows from the
inverse closedness assertion $(d)$, which can be proved using ideas from
\cite{FRS}, where inverse closedness issues of two projections algebras with
coefficients are studied. The proof of assertions $(b)$ and $(c)$ of Theorem
Theorem \ref{nt2} is then standard. \\[3mm]
{\bf Step 4: The index formula.} It remains to prove the index formula $(e)$.
First we have to equip the cylinder $\sT \times \overline{\sR}$ with a suitable
topology, which will be different from the usual product topology. We provide
$\sT$ with the counter-clockwise orientation and $\overline{\sR}$ with the
natural orientation given by the order $<$. Then the desired topology is
determined by the system of neighborhoods $U(t_0, \, \lambda_0)$ of the point
$(t_0, \, \lambda_0) \in \sT \times \overline{\sR}$, defined by
\[
U(t_0, \, -\infty) = \{ (t, \, \lambda) \in \sT \times \overline{\sR} :
|t-t_0| < \delta, t \prec t_0 \} \cup \{(t_0, \, \lambda) \in \sT \times
\overline{\sR} : \lambda < \varepsilon\},
\]
\[
U(t_0, \, +\infty) = \{ (t, \, \lambda) \in \sT \times \overline{\sR} :
|t-t_0| < \delta, t_0 \prec t \} \cup \{(t_0, \, \lambda) \in \sT \times
\overline{\sR} : \varepsilon < \lambda\}
\]
if $\lambda_0 = \pm \infty$ and by
\[
U(t_0, \lambda_0) = \{ (t_0, \, \lambda) \in \sT \times \overline{\sR} :
\lambda_0 - \delta_1 < \lambda < \lambda_0 + \delta_2 \}
\]
if $\lambda_0 \in \sR$, where $\varepsilon \in \sR$ and $\delta, \, \delta_1,
\delta_2$ are sufficiently small positive numbers, and where $t \prec s$
means that $t$ precedes $s$ with respect to the chosen orientation of $\sT$.
Note that the cylinder $\sT \times \overline{\sR}$, provided with the described
topology, is just a homeomorphic image of the cylinder $\sT \times [0, \, 1]$,
provided with the Gohberg-Krupnik topology. The latter has been shown by Gohberg
and Krupnik to be (homeomorphic to) the maximal ideal space of the commutative
Banach algebra $\eT(PC_p)/K(l^p)$; see \cite{GK1} and \cite{BS1}, Proposition
6.28. If one identifies $\sT \times [0, \, 1]$ with $\sT \times \overline{\sR}$,
then the Gelfand transform of a coset $A + K(l^p)$ of $A \in \eT(PC_p)$ is just
the function $\Gamma (A)$ defined in Theorem \ref{t1}.

It is an important point to mention that while the function $\smb A$ for $A \in
\in \eT^0_{k \times k}(PC_p)$ is {\em not} continuous on $\sT \times
\overline{\sR}$ (just consider the south-east entry of $\smb (L(a) P + L(b)
Q)$), the function
\[
(t, \, \lambda) \mapsto  \det \smb A (t, \, \lambda) /(\det a_{22}(t,
\, \infty)  \det a_{22}(t, \, -\infty)
\]
{\em is} continuous on $\sT \times \overline{\sR}$. This non-trivial fact was
observed by Gohberg and Krupnik in a similar situation when studying the
Fredholm theory for singular integral operators with piecewise continuous
coefficients (see \cite{GK3}; an introduction to this topic is also in Chapter
V of \cite{MiP1}).

We will establish the index formula by employing a method which also goes
back to Gohberg and Krupnik and is known as linear extension. This method has
found its first applications in the Fredholm theory of one-dimensional singular
integral equations; see \cite{GoK1,MiP1}. We will use this method in the
slightly different context of Toeplitz plus Hankel operators. Therefore, and
for the readers' convenience, we recall it here.

Let $\cB$ be a unital ring with identity element $e$. With every $h \times
r$-matrix $\beta := (b_{jl})_{j,l = 1}^{h,r}$ with entries in $\cB$, we
associate the element
\begin{equation} \label{ne4}
\el(\beta) = \sum_{j=1}^h b_{j1} \ldots b_{jr} \in \cB
\end{equation}
generated by $\beta$ and call the $b_{jl}$ the generators of $\el(\beta)$.
For each element of this form, there is a canonical matrix $\ext(\beta) \in
\cB_{s \times s}$ with $s = h(r+1) + 1$ with entries in the set $\{0, \, e, \,
b_{jk} : 1 \le j \le h, \, 1 \le k \le r\}$ and with the property that
$\el(\beta)$ is invertible in $\cB$ if and only if $\ext(\beta)$ is invertible
in $\cB_{s \times s}$. Actually, a matrix with this property can be constructed
as follows. Let
\begin{equation} \label{ne5}
\ext(\beta) := \pmatrix{Z & X \cr Y & 0} = \pmatrix{e_{h(r+1)} & 0 \cr W & e}
\pmatrix{e_{h(r+1)} & 0 \cr 0 & \el(\beta)} \pmatrix{Z & X \cr 0 & e}
\end{equation}
where $e_l$ denotes the unit element of $\cB_{l \times l}$,
\[
Z := e_{h(r+1)} + \pmatrix{0 & B_1 & 0 & \cdots & 0 & 0 \cr
0 & 0 & B_2 & \cdots & 0 & 0 \cr
\vdots & \vdots & \ddots & \ddots & \vdots & \vdots \cr
0 & 0 & 0 &  \cdots & 0 & B_r \cr
0 & 0 & 0 & \cdots & 0 & 0}
\]
with $B_j := \diag (b_{1j}, \, b_{2j}, \, \ldots, \, b_{hj})$, $X$ is the
column $- (0, \, \ldots, \, 0, \, e, \, \ldots, \, e)^T$ with $hr$ zeros
followed by $h$ identity elements, $Y$ is the row $(e, \, \ldots, \, e, \,
0, \, \ldots, \, 0)$ with $h$ identity elements followed by $hr$ zeros, and $W
:= (M_0, \, M_1, \, \ldots, \, M_r)$ with $M_0 := (e, \, \ldots, \, e)$
consisting of $h$ identity elements and
\[
M_j := (b_{11} b_{12} \ldots b_{1j}, \, b_{21} b_{22} \ldots b_{2j}, \, \ldots,
\, b_{h1} b_{h2} \ldots b_{hj})
\]
for $j = 1, \, \ldots, \, r$. The matrix $\ext(\beta)$ in (\ref{ne5}) is called
the linear extension of $\el(\beta)$.

Since the outer factors on the right-hand side of (\ref{ne5}) are invertible, it
follows indeed that $\el(\beta)$ is invertible in $\cB$ if and only if its
linear extension $\ext (\beta)$ is invertible in $\cB_{s \times s}$. As a
special case we obtain that if the $b_{jl}$ are bounded linear operators on some
Banach space $B$, then $\el(\beta)$ is a Fredholm operator on $B$ if and only
if $\ext(\beta)$ is a Fredholm operator on $L(B)_{s \times s} = L(B_s)$.
Moreover, $\ind \el(\beta) = \ind \ext(\beta)$ is this case.

We shall apply this observation for $B = l^p(\sZ)_k$ and for the generating
operators
\begin{equation} \label{ne6}
b_{jl} := L(c_{jl}) \, \diag P + L(d_{jl}) \, \diag Q \quad \mbox{with} \;
c_{jl}, \, d_{jl} \in (PC_p)_{k \times k}.
\end{equation}
Put  $\beta := (b_{jl})_{j,l = 1}^{h,r}$, $\gamma := (L(c_{jl}))_{j,l =
1}^{h,r}$ and $\delta := (L(d_{jl}))_{j,l = 1}^{h,r}$. The linear extensions
of $\gamma$ and $\delta$ are Laurent operators again; thus $\ext (\gamma) =
L(c)$ and $\ext (\delta) = L(d)$ with piecewise continuous multipliers $c$ and
$d$. Moreover,
\begin{equation} \label{ne7}
\ext (\beta) = L(c) \, \diag P + L(d) \, \diag Q.
\end{equation}
If $\el (\beta)$ is a Fredholm operator then, by Theorem \ref{nt2} $(a)$, the
matrices $c(t^\pm)$ and $d(t^\pm)$ are invertible for every $t \in \sT$. Hence,
$c$ and $d$ are invertible in $(PC_p)_{ks \times ks}$. This fact together with
the above observation implies that the operator $\el (\beta)$ is Fredholm on
$l^p(\sZ)_k$ if and only if its linear extension $\ext(\beta)$ is Fredholm on
$l^p(\sZ)_{ks}$, which on its hand holds if and only if the Toeplitz operator
$T(d^{-1}c)$ is Fredholm on $l^p_{ks}$, and that the Fredholm indices of the
operators $\el(\beta)$, $\ext(\beta)$ and $T(d^{-1}c)$ coincide in this case.
The symbol of the Toeplitz operator $T(d^{-1}c)$ is the function
\[
\smb (T(d^{-1}c)) (t, \, \lambda) = (d^{-1}c)(t^+) \diag \mu_q(\lambda) +
(d^{-1}c)(t^-) \diag (1-\mu_q(\lambda))
\]
(which stems from the matrix-version of Theorem \ref{t1}), and $\smb
(\ext(\beta)) =: (a_{ij})_{i,j=1}^2$ is related with $\smb (T(d^{-1}c))$ via
\[
\det \smb (T(d^{-1}c)) (t, \, \lambda) = \det (\smb \ext(\beta))(t,
\lambda)/(\det a_{22} (t, \, \infty) \det a_{22} (t, - \infty))
\]
as can be checked directly; see \cite{GoK1,MiP1} for details. This fact can
finally be used to derive the index formula for Fredholm operators of the form
$\el(\beta)$ with the entries of $\beta$ given by (\ref{ne6}). For details we
refer to \cite{GoK1,MiP1} again, where a similar setting is considered.

Since the operators $\el (\beta)$ lie dense in $\eT^0_{k \times k}(PC_p)$, the
index formula for a Fredholm operator in this algebra follows by a standard
approximation argument. To carry out this argument one has to use the
estimate
\[
\| \smb \el(\beta) \| \le M \inf_{K \in K(l^p(\sZ)_k)} \|\el(\beta) + K\|
\]
with $M$ independent of $\beta$, which is an immediate consequence of
(\ref{e271011.2}). \hfill \qed
\section{The index formula for $T+H$-operators}
Our next goal is to provide an index formula for Fredholm operators of the form
$T(a) + H(b)$ on $l^p$ where $a, \, b$ are multipliers in $PC_p$ with a finite
set of discontinuities. We start with a couple of lemmata.
\begin{lemma} \label{l16}
If $a \in C(\sT) \cap M^{\langle p \rangle}$, then $H(a)$ is compact on $l^p$.
\end{lemma}
{\bf Proof.} It is shown in Proposition 2.45 in \cite{BS1} that $C(\sT) \cap
M^{\langle p \rangle} \subseteq C_p$ (in fact it is shown there that the
closure of $C(\sT) \cap M^{\langle p \rangle}$ in the multiplier norm equals
$C_p$) and in Theorem 2.47 that $H(a)$ is compact on $l^p$ if $a \in C_p$.
\hfill \qed \\[3mm]
For a subset $\Omega$ of $\sT$, let $PC(\Omega)$ stand for the set of all
piecewise continuous functions which are continuous on $\sT \setminus \Omega$,
and put $PC_{\langle p \rangle} (\Omega) := PC (\Omega) \cap M^{\langle p
\rangle}$. Thus, $C_{\langle p \rangle} := PC_{\langle p \rangle} (\emptyset) =
C(\sT) \cap M^{\langle p \rangle}$. From 6.27 in \cite{BS1} one concludes that
$PC_{\langle p \rangle} (\Omega) \subseteq PC_p$ if $\Omega$ is finite.

In what follows, we specify $\Omega_0 := \{ \tau_1, \, \ldots, \, \tau_m \}$ to
be a finite subset of $\sT \setminus \{\pm 1\}$ and put $\Omega := \Omega_0 \cup
\{\pm 1\}$. Let $\varphi_0 \in C_{\langle p \rangle}$ be a multiplier which
satisfies $\varphi = \tilde{\varphi}$, takes its values in $[0, \, 1]$, and
is identically 1 on a certain neighborhood of $\{-1, \, 1\}$ and identically 0
on a certain neighborhood of $\Omega_0 \cup \overline{\Omega_0}$. Moreover, we
suppose that $\varphi_0^2 + \varphi_1^2 = 1$ where $\varphi_1 := 1 - \varphi_0$.
\begin{lemma} \label{l17}
Let $c \in PC_{\langle p \rangle} (\{-1, \, 1\})$ and $d \in PC_{\langle p
\rangle} (\Omega_0)$. Then the operators $H(c)T(d) - H(cd \varphi_0)$ and
$T(c)H(d) - H(cd \varphi_1)$ are compact on $l^p$.
\end{lemma}
{\bf Proof.} We write $H(c)T(d) = H(c)T(d) T(\varphi_0) + H(c)T(d) T(\varphi_1)$
with
\begin{eqnarray*}
H(c)T(d) T(\varphi_0) & = & H(c) \left( T(d \varphi_0) - H(d)
H(\widetilde{\varphi_0}) \right) \\
& = & H(cd \varphi_0) - T(c) H(\widetilde{d \varphi_0}) - H(c) H(d)
H(\widetilde{\varphi_0}),
\end{eqnarray*}
\begin{eqnarray*}
H(c)T(d) T(\varphi_1) & = & H(c) T(\varphi_1) T(d) + H(c) \left( T(d)
T(\varphi_1) - T(\varphi_1) T(d) \right) \\
& = & \left( H(c \varphi_1) - T(c) H(\widetilde{\varphi_1}) \right) T(d) \\
& & \qquad + \; H(c) \left(H(d) H(\widetilde{\varphi_1}) -  H(\varphi_1)
H(\tilde{d}) \right).
\end{eqnarray*}
The operators $H(\widetilde{d \varphi_0})$, $H(\widetilde{\varphi_0})$, $H(c
\varphi_1)$, $H(\varphi_1)$ and $H(\widetilde{\varphi_1})$ are compact by Lemma
\ref{l16}, which gives the first assertion. The proof of the second assertion
proceeds similarly. \hfill \qed
\begin{lemma} \label{l18}
Let $a_0, \, b_0 \in PC_{\langle p \rangle} (\{-1, \, 1\})$ and $a_1, \, b_1
\in PC_{\langle p \rangle} (\Omega_0)$. Then the operator
\[
(T(a_0) + H(b_0))(T(a_1) + H(b_1)) - (T(a_0 a_1) + H(a_1 b_0 \varphi_0) +
H(a_0 b_1 \varphi_1))
\]
is compact on $l^p$.
\end{lemma}
{\bf Proof.} We write $(T(a_0) + H(b_0))(T(a_1) + H(b_1))$ as
\begin{eqnarray*}
\lefteqn{T(a_0)T(a_1) + T(a_0)H(b_1) + H(b_0)T(a_1) + H(b_0)H(b_1)} \\
& & = T(a_0 a_1) + K_1 + H(a_0 b_1 \varphi_1) + K_2 + H(b_0 a_1 \varphi_0) +
K_3 + K_4
\end{eqnarray*}
where $K_1 := T(a_0)T(a_1) - T(a_0 a_1)$ and $K_4 := H(b_0)H(b_1) =
T(b_0)T(\widetilde{b_1}) - T(b_0 \widetilde{b_1})$ are compact on $l^p$ by
Proposition 6.29 in \cite{BS1}, and $K_2 := T(a_0)H(b_1) - H(a_0 b_1 \varphi_1)$
and $K_3 := H(b_0)T(a_1) - H(b_0 a_1 \varphi_0)$ are compact by Lemma \ref{l17}.
\hfill \qed \\[3mm]
The following proposition provides us with a key observation; it will allow us
to separate the discontinuities in $\Omega_0$ and $\{-1, \, 1\}$.
\begin{prop} \label{p19}
Let $a, \, b \in PC_{\langle p \rangle} (\Omega)$. If the operator $T(a) +
H(b)$ is Fredholm on $l^p$, then there are functions $a_0, \, b_0 \in
PC_{\langle p \rangle} (\{-1, \, 1\})$ and $a_1, \, b_1 \in PC_{\langle p
\rangle} (\Omega_0)$ such that $T(a_0) + H(b_0)$ and $T(a_1) + H(b_1)$ are
Fredholm operators on $l^p$ and the difference
\[
(T(a_0) + H(b_0))(T(a_1) + H(b_1)) - (T(a) + H(b))
\]
is compact.
\end{prop}
{\bf Proof.} If $T(a) + H(b)$ is Fredholm on $l^p$, then $a$ is invertible in
$PC_p$ by Corollary \ref{c2a} $(a)$. Since the maximal ideal space of $PC_p$ is
independent on $p$ and $a \in PC_{\langle p \rangle}$, one even has $a^{-1} \in
PC_{\langle p \rangle}$.

Let $U$ and $V$ be open neighborhoods of $\{-1, \, 1\}$ and $\Omega_0
\cup \overline{\Omega_0}$, respectively, such that $\clos U \cap \clos V =
\emptyset$. We will assume moreover that $U = U_{-1} \cup U_1$ is the union of
two open arcs such that $\pm 1 \in U_{\pm 1}$, and that $V = V_+ \cup V_-$ is
the union of two open arcs such that $V_+ \subseteq \sT_+^0$ and $V_- \subseteq
\sT \setminus \sT_+^0$. Note that these conditions imply that $\clos U_{-1} \cap
\clos U_1 = \emptyset$.

Now we choose a continuous piecewise (with respect to a finite partition of
$\sT$) linear function $c$ on $\sT$ which is identically 1 on $\clos V$,
coincides with $a$ on $\partial U$, and does not vanish on $\sT \setminus U$.
This function is of bounded total variation; thus $c \in C(\sT) \cap M^{\langle
p \rangle}$, whence $c \in C_p$ as mentioned in the proof of Lemma \ref{l16}.
Put $a_0 := a \chi_U + c \chi_{\sT \setminus U}$. Then $a_0 \in PC_{\langle p
\rangle}$ and $a_0^{-1} \in PC_{\langle p \rangle}$. Further, set $a_1 :=
a_0^{-1} a$. The function $a_1$ is identically 1 on $U$ and coincides with $a$
on $V$. Since $PC_{\langle p \rangle}$ is an algebra, $a_1$ belongs to
$PC_{\langle p \rangle}$. Finally, set $b_0 := b \varphi_0$ and $b_1 := b
\varphi_1$, with $\varphi_0$ and $\varphi_1$ as in front of Lemma \ref{l17}.

The above construction guarantees that $a_0, \, b_0 \in PC_{\langle p \rangle}
(\{-1, \, 1\})$ and $a_1, \, b_1 \in PC_{\langle p \rangle} (\Omega_0)$, and
the operator
\[
(T(a_0) + H(b_0))(T(a_1) + H(b_1)) - (T(a_0 a_1) + H(a_1 b_0 \varphi_0) +
H(a_0 b_1 \varphi_1))
\]
is compact on $l^p$ by Lemma \ref{l18}. The functions $(a_1 - 1) b_0 \varphi_0$
and $(a_0 - 1) b_1 \varphi_1$ vanish identically on a certain neighborhood
of $\Omega$ by their construction. Hence, the Hankel operators $H((a_1 - 1) b_0
\varphi_0)$ and $H((a_0 - 1) b_1 \varphi_1)$ are compact by Lemma \ref{l16},
which implies that the operator
\[
(T(a_0) + H(b_0))(T(a_1) + H(b_1)) - (T(a_0 a_1) + H(b_0 \varphi_0) + H(b_1
\varphi_1))
\]
is compact. Since $a_0 a_1 = a$ and $b_0 \varphi_0 + b_1 \varphi_1 =
b (\varphi_0^2 + \varphi_1^2) = b$, and since $T(a_0) + H(b_0)$ and $T(a_1) +
H(b_1)$ are Fredholm operators on $l^p$ by Theorem \ref{t2}, the assertion
follows. \hfill \qed \\[3mm]
By the previous proposition,
\[
\ind (T(a) + H(b)) = \ind (T(a_0) + H(b_0)) + \ind (T(a_1) + H(b_1)).
\]
Since $H(b_0) \in \eT(PC_p)$ as already mentioned, and since an index formula
for Fredholm operators in $\eT(PC_p)$ is known (see, e.g., 6.40 in \cite{BS1}),
the determination of $\ind (T(a_0) + H(b_0))$ is no serious problem. The
following theorem provides us with a basic step on the way to compute the index
of $T(a_1) + H(b_1)$.
\begin{theo} \label{t20}
Let $a, \, b \in PC_{\langle p \rangle} (\Omega_0)$. If one of the operators
$T(a) \pm H(b)$ is Fredholm on $l^p$, then the other one is Fredholm on $l^p$,
too, and the Fredholm indices of these operators coincide.
\end{theo}
{\bf Proof.} By Corollary \ref{c2a} $(b)$, the operators $T(a) + H(b)$ and $T(a)
- H(b)$ are Fredholm operators on $l^p$ only simultaneously. It remains to prove
that their indices coincide. Recall from the introduction that $T(a) = P L(a) P$
and $H(a) = P L(a) Q J$. Thus, the index equality will follow once we have
constructed a Fredholm operator $D$ such that the difference
\begin{equation} \label{e21}
D(P L(a) P + P L(b) QJ +Q) - (P L(a) P - P L(b) QJ + Q) D
\end{equation}
is compact. The following construction of $D$ is a modification of an idea in
\cite{KaS1}. (Note that the compactness of the operator (\ref{e21}) also
provides an alternate proof of the simultaneous Fredholm property of the
operators $T(a) \pm H(b)$.)

A function $c \in M_p$ is called even (resp. odd) if $c = \tilde{c}$ (resp. $c =
- \tilde{c}$) or, equivalently, if $J L(c) J = L(c)$ (resp. $J L(c) J = -L(c)$).
Every function $c \in C_p$ can be written as a sum of an even and an odd
function in a unique way: $c = (c + \tilde{c})/2 + (c - \tilde{c})/2$. Let
$\theta_o$ and $\theta_e$ be an odd and an even function in $C(\sT) \cap
M^{\langle p \rangle}$, respectively, and assume that $\theta_e$ vanishes at all
points of $\Omega_0$ (and, hence, at all points of $\overline{\Omega_0}$). Put
\begin{equation} \label{e22}
D := P L(\theta_o + \theta_e) P + Q L(\theta_o - \theta_e) Q.
\end{equation}
We will later specify the functions $\theta_o$ and $\theta_e$ such that $D$
becomes a Fredholm operator. First note that
\[
JP L(\theta_o + \theta_e) PJ = - Q L(\theta_o - \theta_e) Q, \quad
JQ L(\theta_o - \theta_e) QJ = - P L(\theta_o + \theta_e) P,
\]
whence $JDJ = -D$ and $JD + DJ = 0$. Next we show that $D$ commutes with the
operator $PL(a)P + P L(b) Q + Q$ up to a compact operator. Since the Toeplitz
operators $PL(\theta_o + \theta_e)P$ and $PL(a)P$ commute modulo a compact
operator, it remains to show that $D$ commutes with $P L(b) Q$ up to
a compact operator. The latter fact follows easily from the identity
\begin{eqnarray*}
\lefteqn{D P L(b) Q - P L(b) Q D} \\
&& = PL(\theta_o + \theta_e) P L(b) Q - P L(b) Q L(\theta_o - \theta_e) Q \\
&& = PL(\theta_o + \theta_e) L(b) Q - PL(\theta_o + \theta_e) Q L(b) Q \\
&& \quad - \; P L(b) L(\theta_o - \theta_e) Q + P L(b) P L(\theta_o - \theta_e)
Q \\
&& = 2 PL(\theta_e b) Q - PL(\theta_o + \theta_e) Q L(b) Q + P L(b) P
L(\theta_o - \theta_e)
Q
\end{eqnarray*}
and from the compactness of the operators $P L(\theta_e b) Q$ and
$PL(\theta_o \pm \theta_e) Q$ by Lemma \ref{l16} (note that $\theta_e b \in
C(\sT) \cap M^{\langle p \rangle}$). The compactness of the operator
(\ref{e21}) is then a consequence of the identity
\begin{eqnarray*}
\lefteqn{D(P L(a) P + P L(b) QJ +Q) - (P L(a) P - P L(b) QJ + Q) D} \\
&& = D P L(a) P - P L(a) P D + D P L(b) QJ + P L(b) QJ D \\
&& = D P L(a) P - P L(a) P D + (D P L(b) Q - P L(b) QD) J
\end{eqnarray*}
and of the compactness of the commutators $[D, PL(a)P]$ and $[D, PL(b)Q]$.

Finally we show that the functions $\theta_e$ and $\theta_o$ can be
specified such that the operator $D$ in (\ref{e22}) is a Fredholm operator on
$l^p$. Set $\hat{\theta}_o (t) := |t^2-1|^2$ for $t \in \sT$. Then
$\hat{\theta}_o$ is an even function in $C^\infty(\sT)$ and $\theta_o :=
\chi_{\sT_+} \hat{\theta}_o - \chi_{\sT_-} \hat{\theta}_o$ is an odd function in
$C(\sT) \cap M^{\langle p \rangle}$. Further,
\[
\theta_e(t) := i \prod_{j=1}^m |t - \tau_j|^2 |t - \overline{\tau_j}|^2, \quad
t \in \sT
\]
defines an even function $\theta_e \in C(\sT) \cap M^{\langle p \rangle}$ which
vanishes at the points of $\Omega_0$. Since $\theta_o$ and $i \theta_e$ are
real-valued functions, we conclude that $\theta_o \pm \theta_e$ are invertible
in $C(\sT) \cap M^{\langle p \rangle}$, which implies that $D$ is a Fredholm
operator as desired. \hfill \qed \\[3mm]
Now we are in a position to derive an index formula for a Fredholm operator of
the form $T(a) + H(b)$ with $a, \, b \in PC_{\langle p \rangle} (\Omega_0)$.
We make use of the well-known identity
\begin{eqnarray} \label{e23}
\lefteqn{\pmatrix{P L(a) P + P L(b) QJ + Q & 0 \cr 0 &  P L(a) P - P L(b) QJ +
Q}} \nonumber \\
&& = \frac{1}{2} \pmatrix{I & J \cr I & -J} \pmatrix{P L(a) P + Q & P L(b) Q
\cr JP L(b) QJ & J(P L(a)P + Q)J} \pmatrix{I & I \cr J & -J},
\end{eqnarray}
where the outer factors in (\ref{e23}) are the inverses of each other.
Thus, if one of the operators $T(a) \pm H(b) = PL(a)P \pm PL(b)QJ$ is a
Fredholm operator, then so is the other, and the Fredholm indices of these
operators coincide. Hence the middle factor
\[
\pmatrix{P L(a) P + Q & P L(b) Q \cr JP L(b) QJ & J(P L(a)P + Q)J} =
\pmatrix{P L(a) P + Q & P L(b) Q \cr Q L(\tilde{b}) P & Q L(\tilde{a}) Q + P}
\]
in (\ref{e23}) is a Fredholm operator, and
\begin{eqnarray*}
\ind(T(a) + H(b)) & = & \frac{1}{2} \, \ind \pmatrix{P L(a) P + Q & P L(b) Q \cr
Q L(\tilde{b}) P & Q L(\tilde{a}) Q + P} \\
& = & \frac{1}{2} \, \ind \pmatrix{P L(a) P & P L(b) Q \cr
Q L(\tilde{b}) P & Q L(\tilde{a}) Q}.
\end{eqnarray*}
For the latter identity note that the operator
\[
A := \pmatrix{P L(a) P + Q & P L(b) Q \cr
Q L(\tilde{b}) P & Q L(\tilde{a}) Q + P} \in L(l^p(\sZ)_2)
\]
has the complementary subspaces $L_1 := \{ (Qx_1, \, Px_2) : (x_1, \, x_2) \in
l^p(\sZ)_2\}$ and $L_2 := \{ (Px_1, \, Qx_2) : (x_1, \, x_2)
\in l^p(\sZ)_2\}$ of $l^p(\sZ)_2$ as invariant subspaces and that $A$ acts on
$L_1$ as the identity operator and on $L_2$ as the operator
\[
A_0 :=  \pmatrix{P L(a) P & P L(b) Q \cr Q L(\tilde{b}) P & Q L(\tilde{a}) Q}.
\]
Let the function $W : \sT \times \overline{\sR} \to \sC$ be defined by
\[
W(t, \, \lambda) = \det \smb A_0 (t, \, \lambda)/(\tilde{a} (t, \,
\infty)\tilde{a}(t, - \infty)).
\]
Since $T(a) + H(b)$ is Fredholm, $W$ does not pass through the origin, and
Theorem \ref{nt2} entails that $\ind A_0 = - \wind W$. Thus,
\[
\ind (T(a) + H(b)) = - \frac{1}{2} \wind W.
\]
We are going to show that actually
\begin{equation} \label{e24}
\ind (T(a) + H(b)) = - \wind_{\sT_+} W,
\end{equation}
where the right-hand side is defined as follows. The compression of $W$ onto
$\sT_+ \times \overline{\sR}$ is a continuous function the values of which form
a closed oriented curve in $\cC$ which starts and ends at $1 \in \sC$ and does
not contain the origin. The winding number of this curve is denoted by
$\wind_{\sT_+} W$. Analogously, we define $\wind_{\sT_-} W$.

For the proof of (\ref{e24}) we suppose for simplicity that $a$ and $b$ have
jumps only at the points $t_1$ and $\overline{t_1}$ where $t_1 \in \sT^0_+$. If
$t$ moves along $\sT_+$ from $1$ to $t_1$ (resp. on $\sT_-$ from $1$ to
$\overline{t_1}$), then the values of $W(t, \, \lambda) = a(t)/\tilde{a}(t) =
a(t)/a(\overline{t})$ move continuously from $1$ to
$a(t_1^-)/a(\overline{t_1}^+)$ (resp. from $1$ to
$a(\overline{t_1}^+)/a(t_1^-)$). Using that $W(t, \, \lambda) = W(\overline{t},
\lambda)^{-1}$ for $t \in \sT \setminus \{-1, \, 1\}$, one easily concludes that
\[
[\arg W]_{1 \to t_1 \subset \sT_+} = [\arg W]_{\overline{t_1} \to 1 \subset
\sT_-}
\]
where the numbers on the left- and right-hand side stand for the increase of
the argument of $W$ if $t$ moves in positive direction along the arc from $1$
to $t_1$ in $\sT_+$ and along the arc from $\overline{t_1}$ to $1$ in $\sT_-$,
respectively. Analogously,
\[
[\arg W]_{-1 \to \overline{t_1} \subset \sT_-} =  [\arg W]_{t_1 \to -1 \subset
\sT_+}.
\]
Consider
\begin{eqnarray*}
\lefteqn{W(t_1, \lambda)/(a(\overline{t_1}^+) a(\overline{t_1}^-))} \\
&& = [a(t_1^+) \mu_q(\lambda) + a(t_1^-)(1 - \mu_q(\lambda))] \,
[a(\overline{t_1}^+) \mu_q(\lambda) + a(\overline{t_1}^-)(1 - \mu_q(\lambda))
] \\
&& \quad - (b(t_1^+) - b(t_1^-)) (b(\overline{t_1}^+) - b(\overline{t_1}^-))
\mu_q(\lambda) (1 - \mu_q(\lambda))
\end{eqnarray*}
and the related expression for $W(\overline{t_1}, \lambda)/(a(t_1^+)
a(t_1^-))$. Then
\[
[\arg W]_{\cC_q(a(t_1^-), \, a(t_1^+))} = [\arg W]_{\cC_q(a(\overline{t_1}^-),
\, a(\overline{t_1}^+))}
\]
because $W(t_1, \, \lambda)/(a(\overline{t_1}^+) a(\overline{t_1}^-)) =
W(\overline{t_1}, \lambda) /(a(t_1^+) a(t_1^-))$. So we arrive at the equality
$\wind_{\sT_+} W = \wind_{\sT_-} W$, whence (\ref{e24}) follows.

Now suppose that $a, \, b \in PC_{\langle p \rangle}$ are continuous on $\sT
\setminus \{-1, \, 1\}$. Then we define a function $W : \sT_+ \times
\overline{\sR}$ by
\[
W(t, \, \lambda) = \left( a(t^+) \mu_q(\lambda) + a(t^-) (1 - \mu_q(\lambda)) +
it (b(t^+) - b(t^-)) \nu_q(\lambda) \right) a^{-1} (\pm 1^\mp)
\]
if $t = \pm 1$ and by $W(t, \, \lambda) = a(t)/a(\overline{t})$ if $t \in
\sT^0_+$. The function $W$ is continuous and determines a closed curve which
starts and ends at $1 \in \sC$. If $T(a) + H(b)$ is a Fredholm operator, then
this curve does not pass through the origin and possesses, thus, a well defined
winding number.

Since $T(a) + H(b)$ is in $\eT(PC_p)$ and the symbol $V : \sT \times
\overline{\sR} \to \sC$ of this operator relative to the algebra $\eT(PC_p)$ is
known (it is just given by
\[
V(t, \, \lambda) = a(t^+) \mu_q(\lambda) + a(t^-) (1 - \mu_q(\lambda)) +
it (b(t^+) - b(t^-)) \nu_q(\lambda)
\]
if $t = \pm 1$ and by $V(t, \, \lambda) = a(t)$ if $t \in \sT \setminus \{-1, \,
1\}$) and since $\ind T(a) = - \wind_{\sT} V$, one can again prove that
$\wind_{\sT} V = \wind _{\sT_+} W$ by comparing the increments of the arguments
as above.

Now we look at the factorization given by Proposition \ref{p19} and denote by
$W_0$ and $W_1$ the above defined function $W : \sT_+ \times \overline{\sR}$ for
the operators $T(a_0) + H(b_0)$ and $T(a_1) + H(b_1)$, respectively. It is easy
to see that $W_0W_1$ coincides with the function $W$ for the operator $T(a) +
H(b)$. Summarizing, we get
\begin{theo}
Let $a, \, b \in PC_{\langle p \rangle}$ and $T(a) + H(b)$ a Fredholm operator
on $l^p$. Then
\[
\ind (T(a) + H(b)) = - \wind_{\sT_+} W_0 - \wind_{\sT_+} W_1 = -\wind_{\sT_+}
W
\]
with $W$, $W_0$ and $W_1$ defined as above.
\end{theo}
\section{The general case} \label{s5}
In this section we want to sketch an approach to derive an index formula for an
arbitrary Fredholm operator $A \in \eT \eH (PC_p)$. With $A$, we associate the
function $W(A) : \sT_+ \times \overline{\sR} \to \sC$ defined by
\[
W(A)(t, \, \lambda) = \left\{
\begin{array}{ll}
\smb A (t, \, \lambda)/\smb A(t, \mp \infty) & \mbox{if} \; t = \pm 1 \\
\det \smb A (t, \, \lambda)/(a_{22} (t, \, \infty) a_{22} (t, \, -\infty)) &
\mbox{if} \; t \neq \pm 1
\end{array}
\right.
\]
where we wrote $\smb A (t, \, \lambda) = (a_{ij} (t, \, \lambda))_{i, j =
1}^2$ for $t \in \sT^0_+$. For $A = T(a) + H(b)$, this definition coincides with
that one from the previous section. 
\begin{theo}
If $A \in \eT \eH(PC_p)$ is a Fredholm operator, then
\begin{equation} \label{e26}
\ind A = - \wind_{\sT_+} W(A).
\end{equation}
\end{theo}
The remainder of this section is devoted to the proof of this theorem. It will
become evident from this proof that $W(A)$ traces out a closed oriented curve
which does not pass through the origin; so the winding number of $W(A)$ is well
defined.

We start with the observation that Theorem \ref{t2} remains true for
matrix-valued multipliers $a, \, b \in (PC_p)_{k \times k}$: just replace
$\mu_q$, $1-\mu_q$ and $\nu_q$ by the corresponding $k \times k$-diagonal
matrices $\diag \mu_q$, $\diag (1-\mu_q)$ and $\diag \nu_q$, respectively. Also
Proposition \ref{pneu1} holds in the matrix setting: If
\[
T(a) + H(b) := (\diag P) L(a) (\diag P) + (\diag P) L(b) (\diag QJ)
\]
is a Fredholm operator, then the identity
\[
\ind (T(a) + H(b)) =  - \wind W(T(a) + H(b))
\]
still holds if one replaces in the above definition of $W$ all scalars by the
determinants of the corresponding matrices. These facts follow in a similar way
as their scalar counterparts.

Now let $a_{jl}, \, b_{jl} \in PC_p$, consider the $h \times r$-matrix $\beta
:= (T(a_{jl}) + H(b_{jl}))_{j,l=1}^{h,r}$, and associate with $\beta$ the
operator
\[
A := \el (\beta) = \sum_{j=1}^h (T(a_{j1}) + H(b_{j1})) \ldots (T(a_{jr}) +
H(b_{jr})) \in \eT \eH (PC_p)
\]
as in (\ref{ne4}). Further set $\gamma := (L(a_{jl}))_{j,l=1}^{h,r}$ and $\delta
:= (L(b_{jl}))_{j,l=1}^{h,r}$. The linear extensions of $\gamma$ and $\delta$
are Laurent operators again; thus $\ext (\gamma) = L(a)$ and $\ext (\delta) =
L(b)$ with certain multipliers $a, \, b \in (PC_p)_{s \times s}$ with $s =
h(r+1) + 1$. Moreover, these extensions are related with the extension of
$\beta$ by
\[
\ext (\beta) = T(\ext (\gamma)) + H(\ext (\delta)) = T(a) + H(b) \in L(l^p_s)
\]
(note that $H(1) = 0$). In Section \ref{s3} we noticed that if $\el (\beta)$ is
Fredholm, then (and only then) $\ext (\beta)$ is Fredholm and $\ind \el (\beta)
=  \ind \ext (\beta)$. Further, if $\el (\beta)$ is a Fredholm operator, then
the matrices $a(t^\pm)$ are invertible for every $t \in \sT$. Hence, $a$ is
invertible in $(PC_p)_{s \times s}$. Now consider
\[
\smb \el (\beta) = \sum_{j=1}^h \smb  (T(a_{j1}) + H(b_{j1})) \ldots \smb
(T(a_{jr}) + H(b_{jr})).
\]
Let $t \neq \pm 1$. Then $\smb (T(a) + H(b))(t, \, \lambda)$ is a matrix of
size $2s \times 2s$. We put the rows and columns of this matrix in
a new matrix according to the following rules: If $j \le h(r+1)+1$, then the $j$
th row of the old matrix becomes the $2j-1$ th row of the new one, whereas if
$j > h(r+1)+1$, the $j$ th row of the old matrix becomes the
$2(j-h(r+1)-1)$ th row of the new matrix. The columns of $\smb (T(a) +
H(b))(t, \, \lambda)$ are re-arranged in the same way. The matrix obtained in
this way is just $\smb \el (\beta) (t, \, \lambda)$. By these manipulations,
\[
\smb \el (\beta) (t, \, \lambda) = \cP \smb (T(a) + H(b)) (t, \, \lambda) \cP^T 
\]
with a certain permutation matrix $\cP$ and its transpose $\cP^T$. Hence, 
\[
\det \smb (T(a) + H(b)) (t, \, \lambda) = \det \smb (\el (\beta)) (t, \,
\lambda)
\]
for $t \neq \pm 1$. For $t = \pm 1$ we do not change the matrix
$\smb (T(a) + H(b))(t, \, \lambda)$.

For $t \neq \pm 1$, we write $\smb (T(a) + H(b)(t, \lambda) = (a_{mn}(t, \,
\lambda))_{m,n=1}^2$ and
\[
\smb (T(a_{jl}) + H(b_{jl}))(t, \, \lambda)) = (a_{mn}^{jl} (t, \,
\lambda))_{m,n=1}^2.
\]
Then
\[
\smb \el(\beta) (t, \, \pm \infty) = \sum_{j=1}^h \prod_{l=1}^r
\pmatrix{a_{11}^{jl} (t, \, \pm \infty) & 0 \cr
0 & a_{22}^{jl} (t, \, \pm \infty)},
\]
and it follows that 
\[
\det a_{22} (t, \, \pm \infty) = \det \ext (\rho(t, \, \pm \infty))
\] 
where $\rho(t, \, \pm \infty) := (a_{22}^{jl} (t, \, \pm \infty))_{j,l =
1}^{hr}$. It is now easy to see that
\[
W(\el (\beta)) (t, \, \lambda) =  W(T(a) + H(b)) (t, \, \lambda) =
W(\ext (\beta)) (t, \, \lambda)
\]
for all $(t, \, \lambda) \in \sT_+ \times \overline{\sR}$, which implies that
$\ind \el (\beta) = - \wind_{\sT_+} W(\el (\beta))$ and, thus, settles the proof
of the index formula (\ref{e26}) for a dense subset of Fredholm operators in
$\eT \eH (PC_p)$.

Finally, we are going to prove estimate (\ref{e241111.1}), i.e., we will show
that there is a constant $M$ such that
\begin{equation} \label{e27}
\|\smb A\|_\infty \le M \inf \{\|A+K\| : K \; \mbox{compact} \}
\end{equation}
for every operator $A \in \eT \eH (PC_p)$. Once this estimate is shown, the
validity of the index formula (\ref{e26}) for an arbitrary Fredholm operator in
$\eT \eH (PC_p)$ will follow by standard approximation arguments as at the end
of Section \ref{s3}.

To prove (\ref{e27}), we consider $\eT \eH (PC_p)$ as a subalgebra of the
smallest closed subalgebra $\eT^0_J (PC_p)$ of $L(l^p(\sZ))$ which contains all
Laurent operators $L(a)$ with $a \in PC_p$, the projection $P$, and the flip
$J$. The homomorphism $\mbox{\rm smb}_t$ defined in Section \ref{s3} cannot be
extended to the algebra $\eT^0_J (PC_p)$ unless $t = \pm 1$. Instead, we are
going to use ideas from \cite{ESi1} and introduce a related family of
homomorphisms $\mbox{smb}_{t, \overline{t}}$ with $t \in \sT^0_ +$
from $\eT^0_J (PC_p)$ onto $(\Sigma^p_1(\sR))_{2 \times 2}$. A crucial
observation (\cite{ESi1}) is that the strong limit 
\begin{equation} \label{e271011.3}
\mbox{smb}_{t, \overline{t}} \, A := \mbox{s-lim}_{n \to \infty}
\pmatrix{A_{t,n,0,0} & A_{t,n,0,1} \cr A_{t,n,1,0} & A_{t,n,1,1}} 
\end{equation}
with $A_{t,n,i,j} := E_n Y_t^{-1} L(\chi_{\sT^+}) J^i A J^j L(\chi_{\sT^+}) Y_t
E_{-n} L_n$ exists for every operator $A \in \eT^0_J (PC_p)$ and every $t \in
\sT_+^0$.
\begin{theo} \label{t021111.2}
Let $t \in \sT^0_+$. Then the mapping $\mbox{\rm smb}_{t, \overline{t}}$ is a
bounded homomorphism from $\eT^0_J (PC_p)$ onto $(\Sigma^p_1(\sR))_{2 \times
2}$. In particular, \\
$(a)$ $\mbox{\rm smb}_{t, \overline{t}} \, P = \diag (\chi_+ I, \, \chi_- I)$
with $\chi_- = 1-\chi_+$, \\
$(b)$ $\mbox{\rm smb}_{t, \overline{t}} \, L(a) = \diag (a(t^+) Q_\sR + a(t^-)
P_\sR, \, a(\overline{t}^-) Q_\sR + a(\overline{t}^+) P_\sR)$ for $a \in
PC_p$, \\
$(c)$ $\mbox{\rm smb}_{t, \overline{t}} \, K = 0$ for every compact operator
$K$, \\
$(d)$ $\mbox{\rm smb}_{t, \overline{t}} \, J = \pmatrix{0 & I \cr I & 0}$.
\end{theo}
{\bf Sketch of the proof.} The existence of the strong limits of the operators
in $(a)$ - $(d)$ and their actual values follow by straightforward computation.
Let us check assertion $(a)$, for instance. For $A=P$, the strong limits of the
diagonal elements of the matrix (\ref{e271011.3}) exist and are equal to $\chi_+
I$ and $\chi_- I$ by Theorem \ref{t271011.1} $(a)$ and since $JPJ = Q$. Now
consider the $01$-entry of that matrix. It is $L(\chi_{\sT^+}) PJ = J
L(\chi_{\sT^-}) Q$ and thus
\begin{eqnarray} \label{e021111.1}
\lefteqn{E_n Y_t^{-1} L(\chi_{\sT^+}) P J L(\chi_{\sT^+}) Y_t E_{-n} L_n}
\nonumber \\
&& = \left(E_n Y_t^{-1} J Y_t E_{-n} \right) \, \left( E_n Y_t^{-1}
L(\chi_{\sT^-}) Q L(\chi_{\sT^+}) Y_t E_{-n} L_n \right). 
\end{eqnarray}
The first factor on the right-hand side is uniformly bounded with respect to
$n$, whereas the second one tends strongly to $0$ by Theorem \ref{t271011.1}
(note that $\chi_{\sT^-} (t) = 0$ for $t \in \sT^0_ +$). Thus, the sequence of
the operators (\ref{e021111.1}) tends strongly to zero. The strong convergence
of the $10$-entry to zero follows analogously.

Another straightforward calculation shows that the mappings
$\mbox{\rm smb}_{t, \overline{t}}$ are algebra homomorphisms and that these
mappings are uniformly bounded with respect to $t \in \sT^0_ +$. Thus, the
mappings $\mbox{\rm smb}_{t, \overline{t}}$ are well-defined on a dense
subalgebra of $\eT^0_J (PC_p)$, and they extend to (uniformly bounded with
respect to $t$) homomorphisms on all of $\eT^0_J (PC_p)$ by continuity.
\hfill \qed \\[3mm]
By assertion $(c)$ of the previous theorem, every mapping $\mbox{\rm smb}_{t,
\overline{t}}$ induces a quotient homomorphism on $\eT^0_J (PC_p)/K(l^p(\sZ))$
in a natural way. We denote this homomorphism by $\mbox{\rm smb}_{t,
\overline{t}}$ again.

Now we are ready for the last step. Let $t \in \sT^0_+$ and $a, \, b \in
PC_p$. From Theorem \ref{t021111.2} we conclude that then the operator
$\mbox{\rm smb}_{t, \overline{t}} (T(a) + H(b))$ is given by the matrix
\[
\pmatrix{\chi_+ (a(t^+) Q_\sR + a(t^-) P_\sR) \chi_+ I &
\chi_+ (b(t^+) Q_\sR + b(t^-) P_\sR) \chi_- I \cr 
\chi_- (b(\overline{t}^-) Q_\sR + b(\overline{t}^+) P_\sR) \chi_+ I &
\chi_- (a(\overline{t}^-) Q_\sR + a(\overline{t}^+) P_\sR) \chi_- I}
\]
acting on $L^p(\sR)_2$. This matrix operator has the complementary subspaces
\[
L_1 := \{ (\chi_- f_1, \, \chi_+ f_2) : f_1, \, f_2 \in L^p(\sR) \}, \;
L_2 := \{ (\chi_+ f_1, \, \chi_- f_2) : f_1, \, f_2 \in L^p(\sR) \} 
\]
of $L^p(\sR)_2$ as invariant subspaces, and it acts as the zero operator on
$L_1$. So we can identify $\mbox{\rm smb}_{t, \overline{t}} (T(a) + H(b))$ with
its restriction to $L_2$, which we denote by $A_0$ for brevity.  

The space $L_2$ can be identified with $L^p(\sR)$ in a natural way. Under this
identification, the operator $A_0$ can be identified with the operator 
\begin{eqnarray*}
\lefteqn{A_1 := \chi_+ (a(t^+) Q_\sR + a(t^-) P_\sR) \chi_+ I + \chi_+ (b(t^+)
Q_\sR + b(t^-) P_\sR) \chi_- I} \\ 
&& + \chi_- (b(\overline{t}^-) Q_\sR + b(\overline{t}^+) P_\sR) \chi_+ I +
\chi_- (a(\overline{t}^-) Q_\sR + a(\overline{t}^+) P_\sR) \chi_- I 
\end{eqnarray*}
which belongs to $\Sigma^p(\sR)$. It is well known (see Section 4.2 in
\cite{RSS1}) and not hard to check that the algebra $\Sigma^p(\sR)$ is
isomorphic to $\Sigma^p_{2 \times 2} (\sR_+)$, where the isomorphism $\eta$ acts
on the generating operators of $\Sigma^p(\sR)$ by
\[
\eta (S_\sR) = \pmatrix{S_{\sR_+} & H_\pi \cr - H_\pi & -S_{\sR_+}} \quad
\mbox{and} \quad 
\eta (\chi_+ I) = \pmatrix{1 & 0 \cr 0 & 0}, 
\]
with $H_\pi$ referring to the Hankel operator 
\[
(H_\pi \varphi)(s) := \frac{1}{\pi i} \int_{\sR_+} \frac{\varphi(t)}{t+s} \,
dt 
\]
on $L^p(\sR_+)$. The entries of the matrix $\eta(A_1)$ are Mellin operators,
and the value of the Mellin symbol of $\eta(A_1)$ at $(t, \, \lambda) \in
\sT_+^0 \times \overline{\sR}$ is the matrix  
\[
\pmatrix{a(t^+) \mu_q(\lambda) + a(t^-) (1-\mu_q(\lambda)) &
(b(t^+) - b(t^-)) \nu_q(\lambda) \cr
(b(\overline{t}^-) - b(\overline{t}^+)) \nu_q(\lambda) &
a(\overline{t}^-) (1-\mu_q(\lambda)) + a(\overline{t}^+) \mu_q(\lambda)},  
\]
which evidently coincides with $\smb (T(a) + H(b)) (t, \, \lambda)$ given in
(\ref{e3}). Summarizing the above arguments we conclude that the homomorphisms 
\[
A + K(l^p) \mapsto (\smb A) (t, \, \lambda)
\]
are uniformly bounded with respect to $(t, \, \lambda) \in \sT_+^0 \times
\overline{\sR}$, which finally implies the estimate (\ref{e27}). \hfill \qed
{\small Authors' addresses: \\[3mm]
Steffen Roch, Technische Universit\"at Darmstadt, Fachbereich
Mathematik, Schlossgartenstrasse 7, 64289 Darmstadt,
Germany. \\
E-mail: roch@mathematik.tu-darmstadt.de \\[3mm]
Bernd Silbermann, Technische Universit\"at Chemnitz, Fakult\"at
Mathematik, 09107 Chemnitz, Germany. \\
E-mail: silbermn@mathematik.tu-chemnitz.de}
\end{document}